\documentstyle[amscd,amssymb,epsf]{amsart}
\newtheorem{thm}{Theorem}[section]
\newtheorem{lem}[thm]{Lemma}
\newtheorem{conj}[thm]{Conjecture}
\newtheorem{prop}[thm]{Proposition}

\theoremstyle{definition}
\newtheorem{defn}[thm]{Definition}

\newtheorem{rem}[thm]{Remark}
\newtheorem{nota}[thm]{Notation}
\newtheorem{conve}[thm]{Convention}

\makeatletter
\def\theequation{\@arabic\c@equation}
\makeatother
\numberwithin{equation}{thm}

\def\intt{\int\limits}

\def\cc{\cal C}
\def\U{{\cal U}}

\def\pt{\Pt}

\def\R{{\Bbb R}}

\def\A{{\cal A}}

\def\V{{\cal V}}
\def\Q{{\Bbb Q}}
\def\Z{{\Bbb Z}}

\def\P{{\Bbb P}}
\def\C{{\Bbb C}}
\def\I{{\cal I}}

\def\aa{{\Bbb A}}

\def\F{{\cal F}}
\def\E{{\cal E}}

\def\oo{\OO}

\def\Ext{{\hbox{\rom{Ext}}}}
\def\Hom{{\hbox{\rom{Hom}}}}
\def\EXT{{\cal Ext}}
\def\HOM{{\cal Hom}}

\def\whsq{\vbox to 5.8pt
{\offinterlineskip\hrule
\hbox to 5.8pt{\vrule height
5.1pt\hss\vrule height 5.1pt}\hrule}}

\def\maps{\longrightarrow}
\def\oo{{\cal O}}

\def\mapr#1{\smash{
 \mathop{\longrightarrow}\limits^{#1}}}

\def\mapl#1{\smash{
 \mathop{\longleftarrow}\limits^{#1}}}

\def\C{{\Bbb C}}
\def\R{{\Bbb R}}
\def\P{{\Bbb P}}

\def\Z{{\Bbb Z}}
\def\Q{{\Bbb Q}}

\def\phi{\varphi}

\def\({\left(}
\def\){\right)}

\def\part{P(n)}

\def\ZZ{{\cal Z}}

\def\pt{{pt}}

\def\<{\langle}
\def\>{\rangle}
\def\Hilb{{\text{\rom{Hilb}}}}

\def\Ext{{\hbox{\rom{Ext}}}}
\def\Sym{{\text{\rm {Sym}}}}

\begin{document}

\title[Donaldson invariants of rational surfaces]
{Wall-crossing formulas, Bott residue formula and the Donaldson invariants
 of rational surfaces}
\author{Geir Ellingsrud}
\address{Mathematical Institute\\University of Oslo\\P.~O.~Box~1053\\
         N--0316 Oslo, Norway}
\email{ellingsr@@math.uio.no}
\author{Lothar G\"ottsche}
\address{Max--Planck--Institut f\"ur Mathematik\\Gottfried--Claren--Stra\ss e
26\\
D-53225 Bonn, Germany}
\email{lothar@@mpim-bonn.mpg.de}

\maketitle\


\section{Introduction}
The  Donaldson invariants
of a smooth $4$ manifold $M$ depend by definition on the choice
of a Riemannian metric.
In case $b^+(M)>1$ they however turn out to be  independent of the metric
 as long as it is generic, and thus they give $C^\infty$-invariants of $M$.
In case $b_+(M)=1$ the invariants have been introduced and studied by
Kotschick in
[Ko]. It turns out that the positive cone of $M$ has a
chamber structure, and Kotschick and Morgan show in \cite{K-M}  that
the invariants only depend on the chamber of the period point of the metric.

Now let $S$ be a smooth algebraic surface with geometric genus
$p_g(S)=0$, irregularity $q(S)=0$, and let $H$ be  an ample divisor on $S$.
Let $M^S_H(c_1,c_2)$ be the moduli space of $H$-Gieseker semistable
rank $2$ sheaves on $S$
with Chern classes $c_1$ and $c_2$.
In the recent paper \cite{E-G} we studied the variation of
$M^S_H(c_1,c_2)$ and that of the corresponding  Donaldson invariants
under change of the ample divisor $H$. For the Donaldson invariants
this corresponds to restricting
our attention from the positive cone of $S$ to the subcone of ample classes.
We imposed a suitable additional condition on the walls between two chambers
and called walls satisfying this condition good walls.

We showed that if the polarisation $H$ passes through a good wall $W$
defined by a cohomology class $\xi\in H^2(S,\Z)$,
then $M^S_H(c_1,c_2)$ changes by a number of flips. Following
\cite{K-M} we wrote the change of the degree $N$ Donaldson invariant
as a sum of contributions  $\delta_{\xi,N}$ with $\xi$ running
through the set of cohomology classes defining $W$. We then used our flip
description to compute $\delta_{\xi,N}$
in terms of
Segre classes of certain standard bundles $\V_{\xi,N}$ over
a Hilbert scheme of points $\Hilb^{d_{\xi,N}}(S\sqcup S)$
on two disjoint copies
of $S$ (here $d_{\xi,N}=(N+3+\xi^2)/4$). We proceeded to compute the leading
 terms of $\delta_{\xi,N}$
explicitely and formulated a conjecture about the precise shape of
$\delta_{\xi,N}$, related to
a conjecture from \cite{K-M}.
We will in future refer to any formula for $\delta_{\xi,N}$
as a wall-crossing formula.

Most of the results of \cite{E-G} were also obtained independently
in \cite{F-Q}, and a flip description of the change of  the
moduli spaces  was obtained independently for varieties of arbitrary
dimension and sheaves of arbitrary rank in \cite{M-W}.
In \cite{H-P} a Feynman path integral aproach  to this problem is developed,
and some of the leading terms of the wall-crossing formulas are determined.

The current paper is a continuation of \cite{E-G}.
We  specialize to the case that the surface $S$ is  rational.
The first advantage is that now almost always all walls are good and so the
formulas from \cite{E-G} almost always apply.

The main reason for restricting our attention to rational surfaces
is that they allow us to use an additional powerful tool:
the Bott residue formula.
A rational surface can always be deformed to a surface
admitting an action of a two-dimensional algebraic torus $\Gamma$
with only a finite number of fixpoints.
As the Donaldson invariants are in particular deformation invariants,
we can assume that $S$ admits such an action of $\Gamma$.
It is easy to see that this action will lift to
the Hilbert schemes $\Hilb^{d_{\xi,N}}(S\sqcup S)$, and that the
standard  bundles
$\V_{\xi,N}$ are equivariant for the induced action. Furthermore also
 the induced action will  only have a finite number of fixpoints, and
 the same is true
for a general $1$-parameter subgroup $T$ of $\Gamma$.
The weights of the action of $T$ on the tangent spaces of
$\Hilb^{d_{\xi,N}}(S\sqcup S)$ and on the fibres of $\V_{\xi,N}$ at the
fixpoints
can be determined explicitely from the corresponding
weights on $S$. So we can apply the Bott
residue formula to this situation and, given $N$ and $\xi$ and the
weights on $S$,
 we  always have an algorithm to compute the change $\delta_{\xi,N}$
explicitely. This algorithm involves very many computations, so
we use  a suitable Maple program.

Now let $S$ be a rational ruled surface with projection
$t:S\maps \P_1$. Let   $F$ be the class of a fibre of $t$
and assume that the intersection number $c_1.F$ is $1$. Then \cite{Q2}
shows
 that, given $c_2\in H^2(S,\Z)$, there always exists a special chamber $\cc_0$
 such that for
$H$ in $\cc_0$ the moduli space
$M^S_H(c_1,c_2)$ is empty. In particular the corresponding
$SO(3)$-invariant is zero on $\cc_0$. This already gives us
an algorithm for computing all the $SO(3)$-invariants corresponding
to  first Chern classes $c_1$ with $c_1.F=1$ on  $S$.
Given a chamber
$\cc$ we obtain the value of the invariant by just summing up all
the changes for all the walls between $\cc$ and $\cc_0$.

At this point we can combine our methods with an additional ingredient:
The blowup formulas, which relate the
Donaldson invariants of an algebraic surface $S$ with those of
the blowup $\widehat S$ of $S$ in a point.
In the case of the projective plane $\P_2$
we obtain an algorithm for computing all the $SO(3)$ and $SU(2)$-invariants.
 Let $\rho:\widehat \P_2\maps \P_2$ be the blow up of $\P_2$
in a point, let $H,F$ and $E$ be the hyperplane class,
the fibre of the projection $\widehat\P_2\maps \P_1$ and  the
exceptional divisor respectively.
We obtain the $SO(3)$-invariants of $\P_2$ by
first computing the  invariants on $\widehat \P_2$ corresponding
to $c_1=\rho^*(H)$ and applying the $SU(2)$-blowup formulas.
Similarly we obtain the
the $SU(2)$-invariants of $\P_2$ by first computing those on $\widehat \P_2$
corresponding to $c_1=E$ and  applying the $SO(3)$-blowup formulas.
Notice that in both cases $c_1.F=1$ on $\widehat \P_2$, so that the
 algorithm of the
previous paragraph applies.
Using a suitable Maple program we have computed all the
$SO(3)$- and $SU(2)$-invariants of $\P_2$ of degree smaller then $50$.

$SO(3)$- and $SU(2)$-invariants of $\P_2$ and rational ruled surfaces
had already been computed by several authors (see e.g. \cite{L-Q}
\cite{E-LP-S} and \cite{K-L}) using a variety of methods.
In \cite{K-L} Kotschick and Lisca have  already made use of the blowup
formulas in combination with the wall-crossing formulas.
Their computations also involve for the first time the $4$-dimensional class.
Their results agree with ours
up to diffenent conventions.
Our paper is partially motivated by and built on  \cite{K-L}.
In particular we found there the correct formulation and the
references for the blowup formulas
in the case $b_+=1$.

We then go back to the wall-crossing  formulas.
Assuming the conjecture from \cite{E-G} about the shape
of $\delta_{\xi,N}$ we are able to determine
(again with a suitable Maple program) the first
$5$ leading terms of $\delta_{\xi,N}$ and using an additional
conjecture even the first $7$ leading terms. Furthermore,
again using the conjecture, we determine $\delta_{\xi,N}$ on
a rational ruled surface for $d_{\xi,N}\le 8$.
By explicitely determining  the corresponding  $\delta_{\xi,N}$
we show that
on a rational ruled surface the conjecture and all the formulas are correct
for all walls $\xi$ and all  $N$ with $N\le 40$
and $d_{\xi,N}\le 8$.

Now we compute the Donaldson invariants for rational ruled surfaces $S$
by again combining the wall-crossing  formulas with the
blowup formulas.  We apply this algorithm to compute all the invariants on
$S$  of degree at most $35$.
The result shows that the special chamber $\cc_0$, where the invariants
corresponding
to first Chern class $c_1$ with $c_1.F=1$ vanish, is also special
for all other $c_1$. We obtain that in the chamber $\cc_0$
the Donaldson invariants  can be expressed as a polynomial
in the linear form  $L_F$ defined by $F$ and the quadratic form $q_S$.
This polynomial is   independent of $S$, and  there is a simple
relationship between the polynomials  for different $c_1$.

Finally we observe that by combining the results obtained so far
with the blowup formulas, we obtain an algorithm for computing
all the $SO(3)$ and $SU(2)$-invariants for all rational
surfaces $S$ for all polarisations in a reasonably big part of the ample cone
of
$S$. This can be seen as a generalization of the result of \cite{K-L}
that  the Donaldson invariants of $\P_2$ and $\P_1\times \P_1$
are determined by the wall-crossing formulas on some blowups.

The explicit computations of the wall-crossing formulas and the
Donaldson invariants of rational surfaces gives us a lot of
empirical data about the shape of these invariants.
We have therefore tried to find some patterns in the results
and so the paper also contains a number of conjectures and questions.
Several of these can already be motivated by the results of \cite{K-L}.

We would like to thank Dieter Kotschick for sending us the preprint
\cite{K-L}, which was quite important for our work,
and also for some useful comments. Furthermore the second author would
like to thank S.A. Str\o mme for a sample Maple program for computations
on Hilbert schemes of points.

\section{Background material}

In this paper let  $S$ be a rational  surface over $\C$.
For such a surface the natural map from the group
of divisors modulo rational equivalence to $H^2(S,\Z)$ is an isomorphism.
So, for $\xi\in H^2(S,\Z)$, we will often write $\oo_S(\xi)$
for the line bundle associated to a divisor with class  $\xi$.

For a polarization  $H$ of $S$ we denote by   $M^S_H(c_1,c_2)$
 the moduli space of
torsion-free sheaves $E$ on $S$ which are H-semistable
(in the sense of Gieseker and Maruyama) of
rank $2$ with $c_1(E)=c_1$ and $c_2(E)=c_2$.

\begin{nota}
For a sheaf $\F$ on a scheme $X$ and a divisor $D$ let $\F(D):=\F\otimes
\oo_X(D)$.
If $X$ is a smooth variety of dimension $n$, we denote the cup product  of two
elements $\alpha$ and $\beta$ in $H^*(X,\Z)$ by $\alpha\cdot\beta$ and the
degree of a class
$\alpha\in H^{2n}(X,\Z)$ by $\int_X\alpha$. For $\alpha,\beta\in H^2(S,\Z)$ let
$\<\alpha\cdot\beta\>:=\int_S\alpha\cdot\beta$. We write
$\alpha^2$ for $\<\alpha\cdot\alpha\>$ and, for $\gamma\in H^2(S,\Z)$, we put
$\<\alpha,\gamma\>:=\<\alpha\cdot \check\gamma\>$, where
$\check\gamma$ is the Poincar\'e dual of $\gamma$.
We denote by $q_S$ the quadratic form on $H_2(S,\Z)$ and, for a
class $\eta\in H^2(S,\Q)$ by $L_\eta$ the corresponding linear form on
$H_2(S,\Q)$.
\end{nota}

\begin{conve} \label{convent}
When we are considering surfaces $S$ and $X$ with a
morphism  $f:X\maps S$, that is either canonical or clear from the context,
then for a cohomology class $\alpha\in H^*(Y,\Z)$ (or a line bundle
$L$ on $Y$) we will very often also denote the pull-back via $f$
by $\alpha$ (resp. $L$). (Very often $f$ will be a sequence of blowups.
In particular if $X$ is a surface which is obtained by $\P_2$ by a number
 of blowups, then we denote by $H$ the pullback of the hyperplane class.
Similarly on $\P_1\times \P_1$ or a variety obtained from $\P_1\times \P_1$
by a number of blowups, we denote by $F$ and $G$ the classes of the fibres
of the projections to the two factors.)
\end{conve}

\bigskip

\subsection{Walls and chambers}
(see \cite{Q1},
\cite{Q2}, \cite{Go}, \cite{K-M} and \cite{E-G}.)

\begin{defn}\label{defwall}
 Let $C_S$ be the ample cone in $H^2(S, \R)$.
For $\xi\in H^2(S,\Z)$ let
$$W^\xi:=C_S\cap\big \{ x\in H^2(S, \R) \bigm| \<x\cdot\xi\>=0\big\}.$$
We shall call $W^\xi$ a wall of type $(c_1,c_2)$, and say that it is
defined by  $\xi$ if the following conditions are satisfied:
\begin{enumerate}
\item  $\xi+c_1$ is divisible by $2$ in
$NS(S)$,
\item $c_1^2-4c_2\le \xi^2<0$,
\item there is a polarisation $H$ with $\< H\cdot \xi\>=0$.
\end{enumerate}

In particular  $d_{\xi,N}:= (4c_2-c_1^2+\xi^2)/4$ is  a nonnegative integer.
An ample divisor $H$ is said to lie in the wall $W$ if $[H]\in W$.
If $D$ is a divisor with $[D]=\xi$, we will also say that
$D$ defines the wall $W$.

A {\it chamber} of type $(c_1,c_2)$ or simply a chamber, is a connected
component of the complement of the union of all the walls of
type
$(c_1,c_2)$.
We will call a wall $W$ {\it good}, if $D+K_S$ is not effective
for any divisor $D$ defining  the wall $W$.
If $(c_1,c_2)$ are given, we call a polarization $L$ of $S$ generic if
it does not lie on a wall of type $(c_1,c_2)$.

On a rational surface $S$ we will call a divisor $L$  {\it good}
if $\<L\cdot K_S\><0$, and we denote by
$C_{S,g}$ the real cone of all good ample divisors.
We see that any wall  $W$  intersecting
$C_{S,g}$ is a good wall.

Let $L_-$ and $L_+$ be two  divisors on $S$.
We denote by
$W_{(c_1,c_2)}(L_-,L_+)$ the set of all $\xi\in H^2(S,\Z)$ defining
a wall
of type $(c_1,c_2)$ and satisfying
$\<\xi\cdot L_-\><0<\<\xi\cdot L_+\>$.
We notice that for $L_-$ and $L_+$ good
all the walls $W^\xi$ defined by $\xi\in W_{(c_1,c_2)}(L_-,L_+)$
are good.

\end{defn}

\subsection
{The change of the Donaldson invariants in terms of Hilbert schemes}

In [Ko] the Donaldson invariants have been introduced for $4$-manifolds
$M$ with $b_+(M)=1$.
In [K-M] it has been shown that in case $b_+(M)=1$, $b_1(M)=0$
they  depend only on the chamber of the
period point of the metric in the positive cone of $H^2(M,\R)$.
We want to use conventions from algebraic geometry,
which differ by a sign from the usual conventions for Donaldson
invariants and furthermore by a factor of a power of $2$ from the conventions
of \cite{Ko}.

\begin{nota}
Let $S$ be a simply connected algebraic surface with $p_g(S)=0$.
Let $N:=4c_2-c_1^2-3$ be a nonnegative integer.
We denote by $A_N(S)$ the set of polynomials of weight $N$
on $H_{2}(S,\Q)\oplus H_{0}(S,\Q)$, where we give weight
$2-i$ to a class in $H_{2i}(S,\Q)$.
Let $\gamma^S_{c_1,N,g}$ be the Donaldson polynomial
of degree $N$ with respect to a generic Riemannian metric $g$
associated to the principal $SO(3)$-bundle $P$ on $S$ whose
second Stiefel-Whitney class
$w_2(P)$ is the reduction of $c_1$ mod $2$ (in the conventions of e.g.
\cite{F-S}).
Then $\gamma^S_{c_1,N,g}$ is a linear map $A_N(S)\maps \Q$.
If $N$ is not congruent to $-c_1^2-3$ modulo $4$, then by definition
$\gamma^S_{c_1,N,g}=0$.
If $g$ is the Fubini-Studi metric associated to
generic ample divisor $L$ on $S$ we denote
$\Phi^{S,L}_{c_1,N}:=(-1)^{(c_1^2+\<c_1\cdot K_S \>)/2}\gamma^S_{c_1,N,g}$.
We denote  $\Phi^{S,L}_{c_1}:=\sum_{N\ge 0} \Phi^{S,L}_{c_1,N}$.
We denote by $\pt\in H_0(S,\Z)$ the class of a point.
Sometimes we will  consider the Donaldson invariants as polynomials on
$H_2(S,\Q)$
by putting
$\Phi^{S,L}_{c_1,N,r}(\alpha):=\Phi^{S,L}_{c_1,N}(\pt^r\alpha^{N-2r})$
for $\alpha\in H_2(S,\Q)$.
\end{nota}

If the moduli space $M_L(c_1,c_2)$  fulfills certain  properties
(in particular there is a universal sheaf $\U$ over
$S\times M_H(c_1,c_2)$),
then for $\alpha_1,\ldots \alpha_r\in H_{2i}(S,\Q)$ we have
$$\Phi^{S,L}_{c_1,N}(\alpha_1\ldots\alpha_r)=\int_{M_L(c_1,c_2)}
\nu(\alpha_1)\cdot\ldots\cdot\nu(\alpha_r)$$
where $\nu(\alpha)=(c_2(\U)-c_1^2(\U)/4)/\alpha$
(\cite{Mo}, \cite{Li}).)

We will use a result from \cite{E-G} (also proved independently in
\cite{F-Q}), We state it only for rational surfaces.
Note that there are some changes in notation.

\begin{defn} \label{wallchange} Let $\xi\in H^2(S,\Z)$ be a
class defining a good wall of type
$(c_1,c_2)$. For $N:=4c_2-c_1^2-3$ we denote
$d_{\xi,N}:=(N+3+\xi^2)/4$, $e_{\xi,N}:=-\<\xi\cdot(\xi-K_S)\>/2+d_{\xi,N}+1$.
Assume now that $(c_1,c_2)$ are fixed.
Let
$$T_\xi:=\Hilb^{d_{\xi,N}}(S\sqcup S)=
\coprod_{n+m=d_{\xi,N}}\Hilb^n(S)\times \Hilb^m(S).$$
be the Hilbert scheme of $d$ points on $2$ disjoint copies of $S$.
Let $q:S\times T_\xi\maps T_\xi$ and $p:S\times T_\xi\maps T_\xi$
be the projections.
Let  $V_\xi$ be the sheaf $p^*(\oo_S(-\xi)\oplus\oo_S(-\xi+K_S))$
on $S\times T_\xi$.
Let $\ZZ_{1}$ (resp.$\ZZ_{2}$) be the subscheme of
$S\times T_\xi$ which restricted to each component
$S\times\Hilb^n(S)\times \Hilb^m(S)$ is the pullback of the universal
subscheme $Z_n(S)$
(resp. $Z_m(S)$)
from the first and second (resp. first and third) factor.
Let $\I_{\ZZ_{1}}$, $\I_{\ZZ_{2}}$
be the corresponding ideal sheaves and $[\ZZ_1]$ and $[\ZZ_2]$
their cohomology classes.
For $\alpha\in H_i(S,\Q)$ let
$\widetilde\alpha:=([\ZZ_1]+[\ZZ_2])/\alpha\in H^{4-i}(T_\xi,\Q)$
Then for $\alpha=\alpha_1\cdot\ldots\cdot\alpha_{N-2r}\pt^r\in A_N(S)$
(with $\alpha_i\in H_2(S,\Q)$) we put
\begin{eqnarray*}
\delta_{\xi,N}(\alpha)&:=&\intt_{T_\xi}\left(\left(\prod_{i=1}^{N-2r}
(\<\alpha_i,\xi/2\>+\widetilde\alpha_i)\right)(-1/4+\widetilde\pt)^r
s(\Ext^1_{q}(\I_{\ZZ_{1}},
\I_{\ZZ_{2}}\otimes V_\xi)\right),\end{eqnarray*}
where $s(\cdot)$ denotes the total Segre class.
We denote $\delta_{\xi}:=\sum_{N>0}\delta_{\xi,N}$.
We will also denote for $\alpha\in H_2(S,\Q)$ by
$\delta_{\xi,N,r}(\alpha):=\delta_{\xi,N}(\pt^r\alpha^{N-2r})$.
\end{defn}

\begin{thm}\label{wallchange1}\cite{E-G},\cite{F-Q}\label{donch1}
Let $S$ be a rational surface.
Let $c_1\in H^2(S,\Z)$ and  $c_2\in \Z$. Let $H_-$ and   $H_+$ be
ample divisors on $S$, such that all the walls defined by elements
of $W_{(c_1,c_2)}(H_-,H_+)$ are good.
Then for all $\alpha\in A_N(S)$  we have
$$\Phi^S_{H_+,N}(\alpha)-\Phi^S_{H_-,N}(\alpha)=
\sum_{\xi\in W_{(c_1,c_2)}(H_-,H_+)}
(-1)^{e_{\xi,N}}
\delta_{\xi}(\alpha).$$
\end{thm}

\subsection{Blowup formulas}

We briefly recall the blowup formulas in the context of algebraic surfaces.
In the case $b_+(S)>1$, when the invariants do not depend on the chamber
structure, they have been shown e.g. in \cite{O}, \cite{L}
and  in the most general form in \cite{F-S}.
In the case $b_+(S)=1$ we cite these results after \cite{K-L}.
By \cite{T} the formulas of \cite{F-S} also hold for $S$ with $b_+(S)=1$,
if the chamber structure is properly taken into account.
Let $S$ be an algebraic surface with $b_+=1$ and let $\epsilon:\widehat S\maps
S$
be the blowup in a point.
Let $E\in H^2(S,\Z)$ be the class of the exceptional divisor.
Let $c_1\in H^2(S,\Z)$ and $c_2\in H^4(S,\Z)$ and put $N=4c_2-c_1^2-3$.
Let $\cc\subset \cc_S$ be a chamber of type $(c_1,c_2)$, let
$\cc_E\subset C_{\widehat S}$ be a chamber of type
$(c_1+E,c_2)$, and  let $ \cc_0\subset C_{\widehat S}$
be a chamber of type $(c_1,c_2)$.
Following \cite{Ko} we say that the  chambers
$\cc$ and $ \cc_E$ (resp. $\cc$ and $ \cc_0$) are related chambers
if
$\epsilon^*(\cc)$ is contained in the closure
$\overline \cc_E$ (resp in $\overline \cc_0$).

\begin{thm}
There are universal polynomials $S_k(x)$ and $B_k(x)$ such that
for all related chambers $\cc$ and $ C_E$  (resp, $\cc$ and $
\cc_0$) as above, all $k\le N$ and all  $\alpha\in A_{N-k}(S)$ we have
\begin{eqnarray}
\label{sobl}\Phi^{\widehat S, \cc_E}_{c_1-E}(\check E^k\alpha)=-
\Phi^{\widehat S, \cc_E}_{c_1+E}(\check E^k\alpha)
&=&\Phi^{S,\cc}_{c_1}(S_k(\pt)\alpha),\\
\label{subl}\Phi^{\widehat S,\cc_0}_{c_1}(\check E^k\alpha)
&=&\Phi^{ S,\cc}_{c_1}(B_k(\pt)\alpha).
\end{eqnarray}
(Note the different sign convention).
The $S_k(x)$ and $B_k(x)$ can  be given in terms of the coefficients
of of the $q$-development of certain $\sigma$-functions.
\end{thm}

We refer to (\ref{sobl}) as $SO(3)$-blowup formulas and to \ref{subl}
as $SU(2)$-blowup formulas.
We will use that the $S_k(x)$ and the $B_k(x)$ are determined by
recursive relations:
(a1) $S_{2k}(x)=0$ for all $k$, (b1) $S_1(x)=1$, $S_3(x)=-x$, $S_5(x)=x^2+2$,
$S_7(x)=-x^3-6x$, (a2) $B_{2k+1}(x) =0$ for all $k$, (b2) $B_{0}(x)=1$,
$B_{2}(x) =0$, $B_{4}(x)=-2$ and, in both cases,  the  recursive relation
\begin{eqnarray*}\label{rec}&&\sum_{i=0}^h{h\choose i}
\big( U_{h+4-i}U_{i}-4U_{h+3-i}U_{i+1}
+6U_{h+2-i}U_{i+2}-4U_{h+1-i}U_{i+3}+U_{h-i}U_{i+4}\big)\\
&&=
-4\sum_{i=0}^h{h\choose i}\big(
xU_{h+2-i}U_{i}+
xU_{h-i}U_{i+2}-2xU_{h+1-i}U_{i+1}+U_{h-i}U_{i}\big),
\end{eqnarray*}
with either $U_i=S_i(x)$ or  $U_i=B_i(x)$
(see e.g. \cite{F-S},\cite{K-L}).

\subsection{The walls for rational surfaces}
Now let $S$ be a rational surface. We want to collect some information
about the set of walls in the ample cone $C_S$.
The following is easy to see:

\begin{rem}
\begin{enumerate}
\item
If $S$ is a rational ruled surface then
$C_S=C_{S,g}$, i.e. all walls are good.
\item
If $S$ is obtained from $\P_2$ by a sequence of blow ups with
exceptional divisors $E_1,\ldots,E_r$ then
$C_{S,g}=C_S\cap\big \{ a(H-a_1E_1-\ldots -a_rE_r) \bigm | a>0,a_i>0,\ \sum_i
a_i<3\big\}$.
\end{enumerate}
\end{rem}

\begin{lem}\label{finwall} For any pair $(H_-,H_+)$ of ample
divisors on a rational surface $S$
and all $c_1\in Pic(S)$ and $c_2\in H^2(S,\Z)$ the set $W_{(c_1,c_2)}(H_-,H_+)$
is finite.
\end{lem}

\begin{pf} The set $\{ tH_-+(1-t)H_+ \ |\ t\in [0,1]\}$ is a compact subset
of $C_S$. Therefore by  \cite{F-M} corollary 1.6  it intersects only
finitely many walls of type $(c_1,c_2)$.
\end{pf}

We now give a list of all walls for $S=\widehat\P_2$ and
$S=\P_1\times\P_1$ which will be used repeatedly in our computations.
We denote by $F=H-E$ the class of a fibre of  $\widehat \P_2\maps \P_1$.
We also denote by $F$ the fibre of the projection to the first factor
of $\P_1\times \P_1$ and by $G$ the class of the fibre of the
projection to the second factor.
The verifications are elementary.

\begin{rem} \label{wallp1p1}
\begin{eqnarray*}
W^{\widehat \P_2}_{0,c_2}(F,H-\delta E)&=&\big\{
2aH-2bE\bigm| b>a>\delta b,\, b^2-a^2\le c_2\big\},\\
W^{\widehat \P_2}_{E,c_2}(F,H-\delta E)&=&\big\{
2aH-(2b-1)E\bigm| b>a>\delta (b-1/2),\,b(b-1)-a^2\le c_2\big\},\\
W^{\widehat \P_2}_{H,c_2}(F,H-\delta E)&=&\big\{
(2a-1)H-2bE\bigm| b\ge a>\delta b+1/2,\,b^2-a(a-1)\le c_2\big\},\\
W^{\widehat \P_2}_{F,c_2}(F,H-\delta E)&=&\big\{
(2a-1)H-(2b-1)E\bigm| b> a>\delta (b-1/2)+1/2,\,b(b-1)-a(a-1)\le c_2\big\},\\
W^{ \P_1\times\P_1}_{0,c_2}(F,F+\delta G)&=&\big\{
2aF-2bG\bigm| 0<b<a\delta,\,2ab\le c_2 \big\},\\
W^{ \P_1\times\P_1}_{F,c_2}(F,F+\delta G)&=&\big\{
(2a-1)F-2bG\bigm| 0<b<(a-1/2)\delta,\,(2a-1)b\le c_2 \big\},\\
W^{ \P_1\times\P_1}_{G,c_2}(F,F+\delta G)&=&\big\{
2aF-(2b-1)G\bigm| 0<b<a\delta+1/2,\,(2b-1)a\le c_2 \big\},\\
W^{ \P_1\times\P_1}_{F+G,c_2}(F,F+\delta G)&=&\big\{
(2a-1)F-(2b-1)G\bigm| 0<b<(a-1/2)\delta+1/2,\,2ab-a-b\le c_2 \big\}.
\end{eqnarray*}
\end{rem}

\subsection{Botts formula}
Now we recall the Bott residue formula (see e.g.
\cite{B},\cite{A-B},\cite{E-S2},\cite{C-L1},\cite{C-L2}).
Let $X$ be a smooth projective variety of dimension $n$
with an algebraic action of the multiplicative group $\C^*$
such that the fixpoint set $F$ is finite.
Differentiation of the action induces a global vector field $\xi\in
H^0(X,T_X)$, and $F$ is precisely the zero locus of $\xi$.
Hence the Koszul complex on the map
$\xi^\vee:\Omega_X\maps \oo_X$ is a locally free resolution of $\oo_F$.
For $i\ge 0$ denote by $B_i$ the cokernel of the Koszul map
$\Omega_X^{i+1}\maps \Omega_X^i$.
It is well known that $H^j(X,\Omega^i_X)=0$ for $i\ne j$.
So there are natural exact sequences for all $i$:
$$0\maps H^i(X,\Omega_X^i)\mapr{p_i} H^i(X,B_i)\mapr{r_i}
H^{i+1}(X,B_{i+1})\maps 0.$$
In particular there are natural maps
$q_i=r_{i-1}\circ\ldots\circ r_0:H^0(F,\oo_F)\maps H^i(X,B_i)$.

\begin{defn}
Let $f:F\maps \C$ be a function and $c\in H^i(X,\Omega^i_X)$.
We say that $f$ represents $c$ if $q_{i+1}(f)=0$ and
$q_i(f)=p_i(c).$
\end{defn}

If $f_1$ represents $a_1\in H^i(X,\Omega^i_X)$ and $f_2$ represents
$a_2\in H^j(X,\Omega^j_X)$, then $f_1f_2$ represents $a_1\cdot a_2
\in H^{i+j}(X,\Omega^{i+j}_X)$.

The following result enables us to compute the degree of
polynomials of weight in the  Chern classes of equivariant vector bundles on
$X$.
  Let $\E$ be an equivariant vector bundle of rank $r$ on $X$.
At each fixpoint $x\in F$ the fibre $\E(x)$ splits as a direct sum of
one-dimensional representations of $\C^*$.
Let $\tau_1(E,x),\ldots\tau_r(\E,x)$ denote the corresponding weights,
and for all $k\ge 0$ let $\sigma_k(\E,x)\in \Z$ be the $k$-th elementary
symmetric function in the $\tau_i(\E,x)$.

\begin{thm}\label{Bottres}
\begin{enumerate}
\item
The $k$-th Chern class $c_k(\E)\in H^k(X,\Omega_X^k)$ of $\E$ can be
represented by
the function $x\mapsto \sigma_k(\E,x)$.
\item
The composition
$H^0(\oo_F)\maps H^n(X,\Omega^n_X)\mapr{res}\C$ maps
$f:F\maps \C$ to $\sum_{x\in F}
f(x)/\sigma_n(T_X,x)$.
\end{enumerate}
\end{thm}

\section{Application of the Bott residue formula}

In this section we want to see how the Bott residue formula
can be used
to compute $\delta_{\xi,N}$ for a class $\xi$ defining a wall
on a rational surface $X$.
Let $\Gamma=\C^*\times \C^*$ be an algebraic 2-torus and let
$\lambda$ and $\mu$ be two independent primitive characters of $\Gamma$.
We identify the representation ring of $\Gamma$ with the ring of Laurent
polynomials in $\lambda$ and $\mu$.
For a variety $Y$ with an action of $\Gamma$ we will denote by
$F_Y$ the set of fixpoints.

\subsection{\label{torS} Actions of a torus on rational surfaces }

We are going to define actions with finitely many fixpoints
of $\Gamma$ on $X=\P_2$, $X=\P_1\times \P_1$
and inductively on surfaces $X=X_r$, where
$X_0=\P_2$ or $X_0=\P_1\times \P_1$ and $X_i$ is the blowup of a
fixpoint of the $\Gamma$-action on $X_{i-1}$.
We also define a lift of the action of $\Gamma$ to all line bundles on
$X$.
These actions  will have the following properties:
\begin{enumerate}
\item Each fixpoint $p\in F_X$ has an invariant neighbourhood $A_p$
isomorphic to $\aa^2=spec([k[x,y])$ on which $\Gamma$ acts by
$t\cdot x=\alpha_px$, $ t\cdot y=\beta_py$ for two independent characters
$\alpha_p$ and $\beta_p$ of $\Gamma$, and the $A_p$ cover $X$.
\item For each line bundle $L\in Pic(X)$ the restriction
$L|_{A_p}$ has a nowhere vanishing section
$s_{L,p}$, with $t\cdot s_{L,p}=\gamma_{L,p}s_{L,p}$ for $\gamma_{L,p}$
a character of $\Gamma$.
\end{enumerate}
{\it (a) The case of $\P_2$.}
Let $T_0,T_1,T_2$ be homogeneous coordinates on $\P_2$.
Let
$\Gamma$ act on  $\P_2$ by $t\cdot T_0= T_0$, $t\cdot T_1=\lambda T_1$,
$t\cdot T_2=\mu T_2$.
The  action of $\Gamma$ has $3$ fixpoints
$p_0:=(1\!\!:\!\! 0\!\! : \!\!0)$, $p_1:=(0\!\! : \!\!1\!\! : \!\!0)$
and $p_2:=(0\!\! : \!\!0\!\! : \!\!1)$.
The  sets  $A_{p_i}:=D(T_i)$ (i.e. the locus  where $T_i\ne 0$)
are affine invariant neighbourhoods.
In appropriate coordinates $x,y$ on $A_{p_0}$
(resp. $A_{p_1}$,$A_{p_2}$), the induced action of $\Gamma$ is
$t\cdot (x,y)=(\lambda x,\mu y)$ (resp. $t\cdot (x,y)=(\lambda^{-1}
x,\mu\lambda^{-1} y)$,
$t\cdot (x,y)=(\mu^{-1}x,\lambda\mu^{-1} y)$).
Furthermore on $A_{p_i}$ the monomial $T_i^n$ defines a
trivializing section of $\oo_{\P_2}(n)$ with
$t\cdot T_0^n=T_0^n$, $t\cdot T_1^n=\lambda^nT_1^n$, $t\cdot T_2^n=\mu^nT_2^n$.

{\it (b) The case of $\P_1\times \P_1$.}
Let $X_0,X_1$ and $Y_0,Y_1$ be homogeneous
 coordinates on the two factors. Let $\Gamma$ act on
$\P_1\times\P_1$ by $t\cdot X_0=X_0$, $t\cdot X_1=\lambda X_1$,
$t\cdot Y_0=Y_0$ and $t\cdot Y_1=\mu Y_1$.
This action has $4$ fixpoints $p_{ij}:=V(X_{1-i})\cap V(Y_{1-j})$
(i.e. the locus  where $X_{1-i}=Y_{1-j}= 0$),
which have affine neightbourhoods $A_{p_{ij}}=D(X_i)\cap D(Y_j)$.
In the appropriate coordinates $x,y$ on $A_{p_{ij}}$ the
action is given by $t\cdot (x,y)=(\lambda^{1-2i}x,\mu^{1-2j}y)$
($i$ and $j\in \{0,1\}$).
Finally a trivializing section of $\oo(n,m)$ on $A_{p_{ij}}$ is
$X_i^n Y_i^m$ with $t\cdot (X_i^nY_j^m)=\lambda^{in}\mu^{jm}X_i^nY_j^m$.

{\it (c) The blowup.}
Now assume that $Y$ is a surface obtained from
$\P_1\times \P_1$ or $\P_2$ by successively blowing up fixpoints
of the action of $\Gamma$, and assume that  the action is extended to $Y$, so
that
it still has finitely many fixpoints, and that the assumptions (1) and (2)
above are satisfied.
Let $p\in F_Y$ be a  fixpoint.
Let $A_{p}$ be an affine neighbourhood of $p$ with coordinates
$x,y$ on which $\Gamma$ acts by $t\cdot (x,y)=(\alpha x,\beta y)$ for
two independent characters $\alpha,\beta$ of $\Gamma$.
Let $X$ be the blowup of $Y$ in $p$, and denote by $E$ the exceptional
divisor and by $\widehat A$ the blow up of $A_{p}$ at $p$.
We can identify $E=\P(\<x,y\>^\vee)$, and the induced action of $\Gamma$ on $E$
has $2$ fixpoints $q_0:=(1\!\! : \!\!0)$ and $q_1:=(0\!\! : \!\!1)$,
which are the fixpoints of $\Gamma$ on $X$ over $p$.
There are  affine neighbourhoods $A_{q_0}=\widehat A\cap D(x)  $
 and $A_{q_1}=\widehat A\cap D(y)$ of $q_0$
and $q_1$ in $X$, with coordinates  $(x,y/x)$ and $(y,x/y)$.
The action
$t\cdot (x,y/x)=(\alpha x,\beta\alpha^{-1}y/x)$,
$t\cdot (y,x/y)=(\beta y,\alpha\beta^{-1}x/y)$ extends the action of
$\Gamma$ on $Y\setminus \{p\}$ to $X$.
Let $L$ be a linebundle on $Y$ with a trivializing section
$s_{L,q}$ near each $q\in F_Y$ with $t\cdot s_{L,q}=\gamma_{L,p_i}s_{L,q}$.
Then $L\otimes \oo(kE)$ has for $i\ne 0$ still $s_{L,q}$ as
a trivializing section near $q$ (with $q\ne p$),
and near $q_0$ (resp. $q_1$)  such a section is
$s_0=s_{L,p}\otimes y^{-k}$ (resp. $s_1=s_{L,p}\otimes x^{-k}$) with
$t\cdot s_0=\gamma_{L,p}\beta^{-k}s_0$ (resp.
 $t\cdot s_1=\gamma_{L,p}\alpha^{-k}s_1$).

\subsection{ The induced action on the Hilbert scheme }

We assume that $S$ is a surface obtained by blowing up $\P_2$ or $\P_1\times
\P_1$ repeatedly, with an action of $\Gamma$ as above.
We fix a positive integer $d$ and want to study the induced action of
$\Gamma $ on the
Hilbert scheme $\Hilb^d(S\sqcup S)$ and on certain "standard bundles"
on $\Hilb^d(S\sqcup S)$, which appear in the wall-crossing formula
 \ref{wallchange}.
The induced action on $\Hilb^d(S\sqcup S)$ is given by
$t\cdot (Y,Z)=(t\cdot Y,t\cdot Z)$, where for a subscheme $Z\subset S$ we
denote by
$t\cdot Z$ the subscheme with ideal
$t\cdot \I_{Z/S}:=\{t\cdot f\ |\ f\in \I_{Z/S}\}$.

Now let $F_S:=\{ p_1,\ldots, p_m\}$ be the set of fixpoints on $S$, and,
for all
$i$, let  $A_i$ be the invariant affine neighbourhood of $p_i$ with coordinates
$x_i,y_i$, such that $t\cdot x_i=\alpha_ix_i,$
$t\cdot y_i=\beta_i y_i$ for two independent
characters $\alpha_i$ and $\beta_i$.
As the characters $\alpha_i$ and $\beta_i$ are independent,
it is easy to see that a subscheme $Z\in \Hilb^n(S)$ is fixed by the
induced action of $\Gamma$ if and only if $supp(Z)\subset F_S$
and if, for all $i$, denoting by $Z_i$ the part of $Z$ with support $p_i$, all
the
ideals
$\I_{Z_i/A_i}$ are generated by monomials in $x_i$ and $y_i$.
We denote by $F_{\Hilb^d(S\sqcup S)}$ the fixpoints
on $\Hilb^d(S\sqcup S)$.

\begin{defn}
A partition of  a nonnegative integer $n$ is a sequence
$\alpha=(a_0,\ldots a_r)$ with
$\alpha_0\ge \ldots \ge a_{r-1}\ge a_r=0$
and $\sum a_i=n$. We identify $(a_0,\ldots a_r)$ and $(a_0,\ldots a_r,0)$.
Let $P_{2m}(d)$ be the set of sequences
$(P_1,\ldots,P_m,Q_1,\ldots Q_m)$ where the $P_i$ and $Q_i$ are all
partitions  of numbers $n_i$ and $m_i$ with $\sum (n_i +m_i)=d$.
We see that $P_{2m}(d)$ and $F_{\Hilb^d(S\sqcup S)}$ are in one-one
correspondence,
with $(P_1,\ldots,P_m,Q_1,\ldots Q_m)$
corresponding to $(Y_1\sqcup\ldots\sqcup Y_m,Z_1\sqcup\ldots\sqcup
Z_m)$, where for $P_i=(a_0,\ldots,a_r)$, $Q_i=(b_0,\ldots,b_r)$
the subschemes $Y_i$ and $Z_i$ are supported at $p_i$ and defined by
$\I_{Y_i/A_i}=(y_i^{a_0},x_iy_i^{a_1}, \ldots x_i^sy_i^{a_s},x^{s+1})$ and
$\I_{Z_i/A_i}=(y_i^{b_0},x_iy_i^{b_1}, \ldots x_i^ry_i^{b_r},x^{r+1})$.
\end{defn}

\subsection{The action on some standard bundles.}
We now want to determine the action of $\Gamma$ on some standard bundles
on $\Hilb^d(S\sqcup S)$ which appear in the wall-crossing formula
\ref{wallchange}.
Let $\xi$ define a good wall.
We denote by $V$ the vector bundle $\oo_S(-\xi)\oplus\oo_S(-\xi+K_S)$.
Then by the results of \cite{E-G} and \cite{F-Q}
$\Ext^1_q(\I_{\ZZ_1},\I_{\ZZ_2}\otimes p^*V)$
is a locally free sheaf on $\Hilb^d(S\sqcup S)$, which
is compatibel with base change, i.e. its fibre over
$(Y,Z)\in \Hilb^d(S\sqcup S)$ is
$\Ext^1(\I_{Y},\I_{Z}\otimes V)$.
Furthermore the $\Gamma$-linearisation of $\oo_S(\xi)$ from \ref{torS}
determines in a canonical way a $\Gamma$-linearisation of
$\Ext^1_q(\I_{\ZZ_1},\I_{\ZZ_2}\otimes p^*V)$.
It also induces an action of $\Gamma$ on $H^1(S,V)$.
Now let $(Y,Z)\in F_{\Hilb^d(S\sqcup S)}$ be a point corresponding to
 $(P_1,\ldots,P_m,Q_1,\ldots Q_m)$.
We will determine the action on the fibre $\Ext^1(\I_{Y},\I_{Z}\otimes V)$.
We denote by $V(p_i)$ the fibre of $V$ over the fixpoint $p_i$
considered as a representation of $\Gamma$.

\begin{lem}\label{stanbott}
For partitions
$P:=(a_0,\ldots,a_r)$, $Q:=(b_0,\ldots,b_r)$ we denote
$$E_{P,Q}(x,y):=\sum_{1\le i\le j\le r}\left(
\sum_{s=a_j}^{a_{j-1}-1} x^{i-j-1}y^{b_{i-1}-s-1}+
\sum_{s=b_j}^{b_{j-1}-1} x^{j-i}y^{s-a_{i-1}}\right).$$
Then in the representation ring of $\Gamma$ we have the identities
\begin{eqnarray}\label{for1} T_{Hilb^d(S\sqcup S)}(Y,Z)&=&
\sum_{i=0}^m  (E_{P_i,P_i}(\alpha_i,\beta_i)+E_{Q_i,Q_i}(\alpha_i,\beta_i)),\\
\label{for2} \Ext^1(\I_{Y},\I_{Z}\otimes V)&=&
H^1(S,V)+\sum_{i=0}^m V(p_i)\cdot E_{P_i,Q_i}(\alpha_i,\beta_i)).
\end{eqnarray}
\end{lem}
\begin{pf}
(\ref{for1}) follows directly from \cite{E-S1}.

\noindent{\it Claim:}
In the representation ring of $\Gamma$ we have the identity
$$\Ext^1(\I_{Y},\I_{Z}\otimes V)=H^1(S,V)+
H^0(S,\EXT^1(\I_Y,\I_Z)\otimes V)+H^0(S,\oo_{Z}\otimes V)
-H^0(S,\HOM(\oo_Y,\oo_Z)\otimes
V).$$
{\it Proof of the Claim:}
As $\xi$ defines a good wall, we have
$H^2(S,\HOM(\I_{Y },\I_{Z})\otimes V)=H^0(S,\HOM(\I_{Z },\I_{Y})
\otimes V^\vee(K_S))=0$ and  $H^0(S,\HOM(\I_{Z },\I_{Y})\otimes V)=0.$
Therefore the low-term exact sequence of the local to global spectral sequence
$H^p(\EXT^q(\I_Y,\I_Z\otimes V))\Rightarrow \Ext^{p+q}(\I_{Y},\I_{Z}\otimes V)$
gives in the representation ring of $\Gamma$
$$\Ext^1(\I_Y,\I_Z(V))
=H^0(S,\EXT^1(\I_Y,\I_Z)\otimes V)+H^1(S,\HOM(\I_Y,\I_Z)\otimes V).$$
We have an exact sequence
$$0\maps \I_{Z}\maps \HOM(\I_{Y},\I_{Z})\maps  \HOM(\oo_{Y},\oo_{Z})\maps
0.$$
So, tensoring by $V$, taking the long exact sequence of cohomology
and using the vanishing of
$H^0(\HOM(\I_Y,\I_Z)\otimes V)$
and  $H^1(\HOM(\oo_{Y },\oo_{Z})\otimes V)$,
we get in the
representation ring of $\Gamma$ the identity
$$H^1(\HOM(\I_{Y},\I_{Z})\otimes V)=H^1(S,\I_Z\otimes V)
-H^0(\HOM(\oo_{Y},\oo_{Z})\otimes V).$$
Finally we use the sequence $0\maps \I_{Z}\otimes V\maps V \maps
\oo_{Z}\otimes V\maps 0$ and  the vanishing of
$H^0(S,V)$ and $H^1(S,\oo_Z\otimes V)$ to replace $H^1(S,\I_Z\otimes V)$
by $H^0(S,\oo_Z\otimes V)+H^1(S,V)$. This shows the claim.

We denote by $\F$ the virtual $\Gamma$-sheaf $\EXT^1(\I_Y,\I_Z)+\oo_{Z}
 -\HOM(\oo_Y,\oo_Z)$. We have to show that $H^0(S,\F\otimes V)
=\sum_{i=0}^m V(p_i)\cdot E_{P_i,Q_i}(\alpha_i,\beta_i).$
If we denote by $\F_i$ the part of $\F$ with support $p_i$, then
$H^0(S,\F\otimes V)=\sum_{i=1}^m H^0(S,\F_i\otimes V).$
We can therefore assume that $supp(Y)=supp(Z)$
is one fixpoint $p$. Let $x$ and $y$ be coordinates near $p$ as before
and $R:=\C[x,y]$.
Let $J:=(y^{a_0},xy^{a_1},\ldots,x^{r+1})$ (resp.
$I:=(y^{b_0},xy^{b_1},\ldots,x^{r+1})$) be the ideal of $Y$ (resp. $Z$).
We denote by $F$ the virtual $R$-$\Gamma$-module corresponding to $\F$.
In the representation ring of $\Gamma$ we have
$$H^0(S,\F\otimes V)=
F\cdot V(p).$$

So we finally have to show that in the representation ring of $\Gamma$ we
have $F=E_{(a_0,\ldots,a_r),(b_0,\ldots b_r)}(\lambda,\mu)$.
The exact sequences
\begin{eqnarray*}
&&0\maps I\maps \Hom_R(J,I)\maps \HOM(R/J,R/I)\maps 0\\
&&0\maps I\maps R\maps R/I\maps 0
\end{eqnarray*}
give
$F=\Ext_R^1(J,I)-\Hom_R(J,I)+R$
in the representation ring of $\Gamma$.

Following \cite{E-S1} we denote by $R[\alpha,\beta]$ the ring $R$
with $\Gamma$-operation defined by $t(x^iy^j):=x^{i-\alpha}y^{j-\beta}$.
We put
$A_0:=\bigoplus_{i=0}^r R[i,a_i],$ $B_0:=\bigoplus_{j=0}^r R[j,b_j],$
$A_1:=\bigoplus_{i=1}^r R[i,a_{i-1}],$ $B_1:=\bigoplus_{j=1}^r R[j,a_{j-1}]$.
Then we have $\Gamma$-equivariant free resolutions
$0\maps A_1\maps A_0\maps J\maps 0$ and $0\maps B_1\maps B_0\maps I\maps 0$.
So the total complex
$$A_0^\vee\otimes B_1\mapr{\alpha}A_1^\vee\otimes B_1\oplus
A_0^\vee\otimes B_0\maps A_1^\vee\otimes B_0$$ associated to the
double complex $\Hom_R(\A_{\bullet},B_{\bullet})$ computes the
$\Ext^i_R(J,I)$, hence
$F=R+A_1^\vee\otimes B_1+A_0^\vee\otimes B_0-A_0^\vee\otimes
B_1-A_1^\vee\otimes B_0.$

Again following \cite{E-S1} we write
$n_i:=(i,a_{i-1})$, $d_i:=(i,a_i)$, $m_j:=(j,b_{j-1})$ and $e_j:=(j,b_j)$.
Then a calculation analogous to \cite{E-S1} shows
$$F=R+\sum_{{1\le i\le r}\atop {0\le j\le r}}
R[e_j-n_i]-\sum_{{1\le i\le r}\atop {1\le j\le r}} R[m_j-n_i]
-\sum_{{0\le i\le r}\atop {0\le j\le r}}R[e_j-d_i]+
\sum_{{0\le i\le r}\atop {1\le j\le r}}R[m_j-d_i].$$
Putting
\begin{eqnarray*}
K_{i,j}&:=&R[m_j-d_{i-1}]-R[m_j-n_i]-R[e_j-d_{i-1}]+R[e_j-n_i],\\
L_{i,j}&:=&R[m_i-d_{j}]-R[m_i-n_j]-R[e_{i-1}-d_{j}]+R[e_{i-1}-n_j],
\end{eqnarray*}
a calculation analogous to \cite{E-S1} gives
$$F=\sum_{1\le i\le j\le r}
(K_{i,j}+L_{i,j}),\
K_{i,j}=\sum_{s={a_i}}^{a_{i-1}} \lambda^{i-j-1}\mu^{b_{i-1}-s-1}\hbox{ and }
L_{i,j}=\sum_{s={b_j}}^{b_{j-1}} \lambda^{j-i}\mu^{s-a_{i-1}},
$$
and the result follows.
\end{pf}

We want to use the easy fact that representation of cohomology classes
is  compatible with equivariant pullback:
 Let  $X$ and $Y$ be smooth projective varieties with an action
 of $\C^*$ with finitely many fixpoints
and let $\mu:X\maps Y$ be an equivariant surjective  morphism.
Then $\mu$ induces a morphism $\mu|_{F_X}:F_X\maps F_Y$.

\begin{lem}
\label{compat}
$f\in \oo_{F_Y}$ represents a cohomology class
$c\in H^j(Y,\Omega_Y^j)$ if and only if
$(\mu|_{F_X})^*f$ represents $\mu^*c$.
\end{lem}

\begin{lem}\label{altild} Let $\alpha\in H^k(S,\Z)$ be a class represented by
$f:F_S\maps \C$. Then  $\widetilde \alpha$ (see \ref{wallchange})
on $\Hilb^d(S\sqcup S)$
is represented by
$$\widetilde f:P_{2m}(d)\maps \C,
((P_i),(Q_i))\mapsto \prod_{i=1}^m (n_i+m_i) f(p_i),$$
where $P_i\in P(n_i)$ and $Q_i\in P(m_i)$.
\end{lem}
\begin{pf}
Let $$\Hilb^{d-1,d}(S\sqcup S):=\big\{(Z_{d-1},Z_d)\in \Hilb^{d-1}(S\sqcup S)
\times \Hilb^{d}(S\sqcup S)\bigm | Z_{d-1}\subset Z_d\big\}$$
with the reduced induced structure. Then
$\Hilb^{d-1,d}(S\sqcup S)$ is smooth and we have a diagram
$$\Hilb^d(S\sqcup S)\mapl{\phi}\Hilb^{d-1,d}(S\sqcup S)\mapr{\psi}
(S\sqcup S)\times \Hilb^{d-1}(S\sqcup S)\mapr{\eta} S\times \Hilb^{d-1}(S\sqcup
S)$$
Here $\psi$ is the blowup along the universal family $Z_d(S\sqcup S)$
\cite{E}
and $\eta$
is induced by the identity map on $S$ and $\Hilb^d(S\sqcup S)$.
It is easy to see from the definitions that
$\phi^*\widetilde\alpha=\psi^*\eta^*(p_1^*\alpha+p_2^*\bar \alpha)$,
where $p_1$ and $p_2$ are the projections of $S\times \Hilb^{d-1}(S\sqcup S)$
onto its two factors and $\bar \alpha$ is the class corresponding
to $\widetilde \alpha$ if we replace $d$ by $d-1$.
It is easy to see that $\phi$, $\psi$ and $\eta$ are equivariant for the
natural lifts of the action of $\Gamma$ on $S$, furthermore
the fixpoint sets $F_{\Hilb^d(S\sqcup S)}$, $F_{\Hilb^{d-1,d}(S\sqcup S)}$
and $F_{S\times \Hilb^{d-1}(S\sqcup S)}$ are all finite.
In fact
 we can identify
$$F_{\Hilb^{d-1,d}(S\sqcup S)}
=\big \{((S_i,T_i),(P_i,Q_i))\in P_{2m}(d-1)\times P_{2m}(d)\bigm|
P_i\ge S_i, \ Q_i\ge T_i \hbox{ for all } i\big\},$$
where for partitions $P=(a_1,\ldots ,a_r),$ $Q=(b_1,\ldots b_r)$ we denote by
$P\ge Q$ that $a_i\ge b_i$ for all $i$.
Obviously
$F_{S\times \Hilb^{d-1}(S\sqcup S)}=
F_S\times P_{2m}(d-1)$ and with this identification
$\phi$ and $\eta\circ\psi$ are the obvious maps.
Now, applying lemma \ref{compat} to $\phi$ and $\eta\circ\psi$,
 the result follows by easy induction.
\end{pf}

We can now put our results together:

\begin{nota} Fix  a one-parameter subgroup $T$ of $\Gamma$.
Let $\xi$ define a good wall on $S$. For any line bundle $L$ on $S$
denote by $w_i(L)$ the weight of  the induced action of $T$ on
the fibre $L(p_i)$.
Let $L_1$ and $L_2$ be two line bundles with
$\<L_1\cdot L_2\>=\pt$ (e.g. if $S$ is a blow up of $\P_2$ then we take
$L_1=L_2=H$).
Furthermore denote by $w(x_i)$, $w(y_i)$ the weight of
the action of $T$ on $x_i$, $y_i$.
We denote for partitions $P=(a_0,\ldots a_r)$ and $Q=(b_0,\ldots b_r)$
 of numbers $n$ and $m$
\begin{eqnarray*}
\overline F_{P,Q}(u,v)&:=&\prod_{1\le i\le j \le r}
\prod_{s=a_j}^{a_{j-1}-1} ((i-j-1)u+(b_{i-1}-s-1)v)
\prod_{s=b_j}^{b_{j-1}-1} ((j-i)u+(s-a_{i-1})u)\\
F^z_{P,Q}(u,v,t)&:=&\prod_{1\le i\le j \le r}
\prod_{s=a_j}^{a_{j-1}-1} (1+z((i-j-1)u+(b_{i-1}-s-1)v+t))\\&&\qquad
\prod_{s=b_j}^{b_{j-1}-1} (1+z((j-i)u+(s-a_{i-1})v+t))
\end{eqnarray*}
By lemma \ref{stanbott}, when putting the correct weights
$F_{P,Q}(u,v)$ will represent the top Chern class of $\Hilb^d(S\sqcup S)$
and $F^z_{P,Q}(u,v,t)$ the total Chern class of
$\Ext^1_q(\I_{\ZZ_1},\I_{\ZZ_2}\otimes p^*V)$.
\end{nota}

\begin{thm}\label{Botthilb}
Let $\alpha_1,\ldots,\alpha_{N-2r}\in H^2(S,\Z)$.
 If $T$ is sufficiently general, then
\begin{eqnarray*}&&\label{bottformel}
\delta_{\xi}(\alpha_1\alpha_2\ldots\alpha_{N-2r}\pt^r)=
{\hbox{\rm Coeff}}_{z^{2d}}\Bigg(\sum_{((P_i),(Q_i))\in P_{2m}(d)}\\&&
\Bigg(
\prod_{k=1}^{N-2r}\Big(\<\xi,\alpha_k\>/2+\sum_{i=1}^m
w_i(\alpha_k)(n_i+m_i)z\Big)
\Big(-1/4+\sum_{i=1}^m w_i(L_1)w_i(L_2)(n_i+m_i)z^2\Big)^r\cdot\\
&&\Bigg(\prod_{i=1}^m\Big(\overline F_{(P_i,P_i)}(w(x_i),w(y_i))
\overline F_{Q_i,Q_i}(w(x_i),w(y_i))\cdot\\
&&F^z_{P_i,Q_i}(w(x_i),w(y_i),-w_i(\xi))
F^z_{P_i,Q_i}(w(x_i),w(y_i),-w_i(\xi))+w_i(K_S))\Bigg)^{-1}\Bigg).
\end{eqnarray*}
\end{thm}

\begin{pf}
The Chern classes of $\V_\xi=\Ext^1_q(\I_{\ZZ_1},\I_{\ZZ_2}\otimes p^*V)$
are the same as those of the virtual bundle
$\V_\xi-H^1(S,V)\otimes \oo_{T_{\xi}}$.
Therefore the result just follows by putting together
lemma \ref{stanbott}, lemma \ref{altild} and applying the
Bott residue formula \ref{Bottres}. Notice that $T$ is sufficiently
general if none of the denominators  vanish.
\end{pf}

This formula can be implemented as a Maple program.

\section{The Donaldson invariants of the
projective plane}

In this section we want to compute the $SU(2)$- and the $ SO(3)$-invariants of
the projective plane $\P_2$ by first computing on the blowup $\widehat \P_2$
and then using the blowup formulas.

In order to get started we need the following easy result of \cite{Q2}:

\begin{lem}\label{vancham}
Let $S$ be a rational ruled surface, $F$ the class of a fibre and $E$ the class
of a section. Fix $(c_1,c_2)\in H^2(S,\Z)\times H^4(S,\Z)$
with $\<c_1\cdot F\>=1$. Then, for all $\epsilon>0$
which are sufficiently small,
we have $M_{F+\epsilon E}(c_1,c_2)=\emptyset$. In particular we get
for $N:=4c_2-c_1^2-3$ that $\Phi^{S,F+\epsilon E}_{c_1,N}=0$.
\end{lem}

We will denote by $E$ the exceptional divisor on $\widehat \P_2$
and by $H$ the (pullback of) the hyperplane class on $\P_2$.

\subsection{The $SU(2)$-case}
We first consider the $SO(3)$-invariants on $\widehat \P_2$ with respect to
Chern classes
$(E,c_2)$ and put $N:=4c_2-3$.
For $0<\epsilon <<1$ the polarisation $L_\epsilon:=H-\epsilon E$ of $\widehat
\P_2$ lies in a chamber of type $(E,c_2)$ which is related to the polarisation
$H$ of $\P_2$. Thus  (\ref{sobl}) gives
$$\Phi^{\P_2,H}_{0,N}(H^{N-2r}\pt^r)=
\Phi^{\widehat\P_2,L_{\epsilon}}_{E,N+1}(\check E\check H^{N-2r}\pt^r).$$
On the other hand we know by lemma \ref{vancham} that
$\Phi^{\widehat \P_2,L_{1-\epsilon}}_{E,N+1}=0$, for
$L_{1-\epsilon}:=H-(1-\epsilon) E$.
Thus we get
$$\Phi^{\P_2,H}_{0,N}(H^{N-2r}\pt^r)=
\sum_{\xi\in W^{\widehat \P_2}_{E,c_2}(H-E,H)}(-1)^{e_{\xi,N+1}}
\delta_{\xi,N+1}(\check E \check H^{N-2r}\pt^r),$$
where $W^{\widehat \P_2}_{E,c_2}(H-E,H)$ is known by remark \ref{wallp1p1}.
Now we compute the $\delta_{\xi,N+1}(\check E \check H^{N-2r}pt^r)$
 with a maple program
using the Bott residue theorem (i.e. theorem \ref{Botthilb}).
For $N:=4i+1$ we denote
$$A_N:=\sum_{j=0}^{2N} \Phi^{\P_2,H}_{0,N}(\check H^{N-2j}\pt^j)h^{N-2j}p^j.$$
Then our result is:

\begin{thm}\label{p2su}
The $SU(2)$-invariants of $\P_2$ are
\par
{\tolerance 10000\small\parindent 0pt\raggedright
$A_1=-{\frac {3\,h}{2}}$,
$A_5={h}^{5}-p{h}^{3}-{\frac {13\,{p}^{2}h}{8}}$,
$A_9=3\,{h}^{9}+{\frac {15\,p{h}^{7}}{4}}
-{\frac {11\,{p}^{2}{h}^{5}}{16}}-
{\frac {141\,{p}^{3}{h}^{3}}{64}}-{\frac {879\,{p}^{4}h}{256}}$,
\par
$A_{13}=54\,{h}^{13}+24\,p{h}^{11}+{\frac {159\,{p}^{2}{h}^{9}}{8}}+{\frac {51
\,{p}^{3}{h}^{7}}{16}}-{\frac {459\,{p}^{4}{h}^{5}}{128}}-{\frac {1515
\,{p}^{5}{h}^{3}}{256}}-{\frac {36675\,{p}^{6}h}{4096}}$,
\par
$A_{17}=2540\,{h}^{17}+694\,p{h}^{15}+{\frac {487\,{p}^{2}{h}^{13}}{2}}+{
\frac {2251\,{p}^{3}{h}^{11}}{16}}+{\frac {2711\,{p}^{4}{h}^{9}}{64}}-
{\frac {5\,{p}^{5}{h}^{7}}{16}}-{\frac {3355\,{p}^{6}{h}^{5}}{256}}-{
\frac {143725\,{p}^{7}{h}^{3}}{8192}}-{\frac {850265\,{p}^{8}h}{32768}
}$,
\par
$
A_{21}=233208\,{h}^{21}+45912\,p{h}^{19}+10625\,{p}^{2}{h}^{17}+3036\,{p}^{3}
{h}^{15}+{\frac {41103\,{p}^{4}{h}^{13}}{32}}+{\frac {1741\,{p}^{5}{h}
^{11}}{4}}+{\frac {5619\,{p}^{6}{h}^{9}}{64}}
-{\frac {20379\,{p}^{7}{h
}^{7}}{1024}}-{\frac {754141\,{p}^{8}{h}^{5}}{16384}}-{\frac {904239\,
{p}^{9}{h}^{3}}{16384}}-{\frac {10504593\,{p}^{10}h}{131072}}$,
\par
$A_{25}=35825553\,{h}^{25}+{\frac {21975543\,p{h}^{23}}{4}}+{\frac {15224337\,
{p}^{2}{h}^{21}}{16}}+{\frac {12159687\,{p}^{3}{h}^{19}}{64}}+{\frac {
11618625\,{p}^{4}{h}^{17}}{256}}+{\frac {15077511\,{p}^{5}{h}^{15}}{
1024}}
+{\frac {19602561\,{p}^{6}{h}^{13}}{4096}}+{\frac {20676279\,{p}
^{7}{h}^{11}}{16384}}+{\frac {11107665\,{p}^{8}{h}^{9}}{65536}}-{
\frac {28437201\,{p}^{9}{h}^{7}}{262144}}-{\frac {169509159\,{p}^{10}{
h}^{5}}{1048576}}-{\frac {757633329\,{p}^{11}{h}^{3}}{4194304}} -{
\frac {4334081031\,{p}^{12}h}{16777216}}$,
\par
$A_{29}=8365418914\,{h}^{29}+1047342410\,p{h}^{27}+{\frac {1157569571\,{p}^{2}
{h}^{25}}{8}}+{\frac {357034013\,{p}^{3}{h}^{23}}{16}}+{\frac {
499796309\,{p}^{4}{h}^{21}}{128}}+{\frac {25506259\,{p}^{5}{h}^{19}}{
32}}+{\frac {423516455\,{p}^{6}{h}^{17}}{2048}}+{\frac {245576651\,{p}
^{7}{h}^{15}}{4096}+{\frac {537423737\,{p}^{8}{h}^{13}}{32768}}+{
\frac {118590907\,{p}^{9}{h}^{11}}{32768}}}+{\frac {131266019\,{p}^{10}
{h}^{9}}{524288}}-{\frac {498648655\,{p}^{11}{h}^{7}}{1048576}}-{
\frac {4800905323\,{p}^{12}{h}^{5}}{8388608}}-{\frac {2551074181\,{p}^
{13}{h}^{3}}{4194304}}-{\frac {115237180987\,{p}^{14}h}{134217728}}$,
\par
$A_{33}=2780195996868\,{h}^{33}+293334321858\,p{h}^{31}+{\frac {67261095005\,{
p}^{2}{h}^{29}}{2}}+{\frac {67539891519\,{p}^{3}{h}^{27}}{16}}+{\frac
{37480404303\,{p}^{4}{h}^{25}}{64}}+{\frac {1455758501\,{p}^{5}{h}^{23
}}{16}}+{\frac {258640401\,{p}^{6}{h}^{21}}{16}}+{\frac {14239101477\,
{p}^{7}{h}^{19}}{4096}+{\frac {14274421501\,{p}^{8}{h}^{17}}{16384}}+
{\frac {7420640919\,{p}^{9}{h}^{15}}{32768}}}+{\frac {7179481275\,{p}^{
10}{h}^{13}}{131072}}+{\frac {10830752675\,{p}^{11}{h}^{11}}{1048576}}
-{\frac {218792349\,{p}^{12}{h}^{9}}{4194304}}-{\frac {8125524657\,{p}
^{13}{h}^{7}}{4194304}}-{\frac {34310453897\,{p}^{14}{h}^{5}}{16777216
}}-{\frac {561608698989\,{p}^{15}{h}^{3}}{268435456}}-{\frac {
3135392459541\,{p}^{16}h}{1073741824}}$,
\par
$A_{37}=1253558847090600\,{h}^{37}+114049802084088\,p{h}^{35}+11151310348527\,
{p}^{2}{h}^{33}+{\frac {2356053433779\,{p}^{3}{h}^{31}}{2}}+{\frac {
4330247481231\,{p}^{4}{h}^{29}}{32}}+{\frac {136302640305\,{p}^{5}{h}^
{27}}{8}}+{\frac {19010805303\,{p}^{6}{h}^{25}}{8}}+{\frac {5969251539
\,{p}^{7}{h}^{23}}{16}}+{\frac {560820990153\,{p}^{8}{h}^{21}}{8192}}+
{\frac {120608801397\,{p}^{9}{h}^{19}}{8192}}+{\frac {227818311585\,{p
}^{10}{h}^{17}}{65536}}+{\frac {108448319625\,{p}^{11}{h}^{15}}{131072
}}+{\frac {380576939595\,{p}^{12}{h}^{13}}{2097152}}+{\frac {
30600598425\,{p}^{13}{h}^{11}}{1048576}}-{\frac {24269489295\,{p}^{14}
{h}^{9}}{8388608}}-{\frac {128245327215\,{p}^{15}{h}^{7}}{16777216}}-{
\frac {3958008534873\,{p}^{16}{h}^{5}}{536870912}}-{\frac {
3928400321367\,{p}^{17}{h}^{3}}{536870912}}-{\frac {43427017514031\,{p
}^{18}h}{4294967296}}$,
\par
$A_{41}=739328941273642584\,{h}^{41}+59025071651407086\,p{h}^{39}+{\frac {
10043994491097141\,{p}^{2}{h}^{37}}{2}}+{\frac {3657609707934747\,{p}^
{3}{h}^{35}}{8}}+{\frac {1431875545375857\,{p}^{4}{h}^{33}}{32}}+{
\frac {605936823728685\,{p}^{5}{h}^{31}}{128}}+{\frac {279056393051655
\,{p}^{6}{h}^{29}}{512}}+{\frac {141081414974709\,{p}^{7}{h}^{27}}{
2048}}+{\frac {79278076181247\,{p}^{8}{h}^{25}}{8192}}+{\frac {
50786721057147\,{p}^{9}{h}^{23}}{32768}}+{\frac {37720541142561\,{p}^{
10}{h}^{21}}{131072}}+{\frac {31629775641351\,{p}^{11}{h}^{19}}{524288
}}+{\frac {28182330795381\,{p}^{12}{h}^{17}}{2097152}}+{\frac {
24970918823121\,{p}^{13}{h}^{15}}{8388608}}+{\frac {20134588411731\,{p
}^{14}{h}^{13}}{33554432}}+{\frac {10829868897921\,{p}^{15}{h}^{11}}{
134217728}}-{\frac {9966987184029\,{p}^{16}{h}^{9}}{536870912}}-{
\frac {63772382356485\,{p}^{17}{h}^{7}}{2147483648}}-{\frac {
230109186457887\,{p}^{18}{h}^{5}}{8589934592}}-{\frac {891115248823257
\,{p}^{19}{h}^{3}}{34359738368}}-{\frac {4881669867807723\,{p}^{20}h}{
137438953472}}$,
\par
$A_{45}=554194295294679879984\,{h}^{45}+39362965900726633056\,p{h}^{43}+
2959973227900487391\,{p}^{2}{h}^{41}+{\frac {472790595509415927\,{p}^{
3}{h}^{39}}{2}}+{\frac {321971749339669677\,{p}^{4}{h}^{37}}{16}}+{
\frac {29333848336377675\,{p}^{5}{h}^{35}}{16}}+{\frac {
45985006658693745\,{p}^{6}{h}^{33}}{256}}+{\frac {9746506854402795\,{p
}^{7}{h}^{31}}{512}}+{\frac {9002198684193567\,{p}^{8}{h}^{29}}{4096}}
+{\frac {285870961148823\,{p}^{9}{h}^{27}}{1024}}+{\frac {
2602166222135403\,{p}^{10}{h}^{25}}{65536}}+{\frac {844639104181119\,{
p}^{11}{h}^{23}}{131072}}+{\frac {1247528861178057\,{p}^{12}{h}^{21}}{
1048576}}+{\frac {252908840631825\,{p}^{13}{h}^{19}}{1048576}}+{\frac
{855800427325917\,{p}^{14}{h}^{17}}{16777216}}+{\frac {356809084455699
\,{p}^{15}{h}^{15}}{33554432}}+{\frac {532554096813723\,{p}^{16}{h}^{
13}}{268435456}}+{\frac {28716835828749\,{p}^{17}{h}^{11}}{134217728}}
-{\frac {398090966021613\,{p}^{18}{h}^{9}}{4294967296}}-{\frac {
983201250012705\,{p}^{19}{h}^{7}}{8589934592}}-{\frac {
6735895639287969\,{p}^{20}{h}^{5}}{68719476736}}-{\frac {
51089630811025563\,{p}^{21}{h}^{3}}{549755813888}}-{\frac {
1110523927325938473\,{p}^{22}h}{8796093022208}}
$,
\par
$A_{49}=
515844680321852815028832\,{h}^{49}+32961783591325975299120\,p{h}^{47}+
2217961739175425919036\,{p}^{2}{h}^{45}+$\par${\frac {315133878583666626003
\,{p}^{3}{h}^{43}}{2}}+{\frac {94814551550640933927\,{p}^{4}{h}^{41}}{
8}}+{\frac {3787623211889821185\,{p}^{5}{h}^{39}}{4}}+{\frac {
645212007742318605\,{p}^{6}{h}^{37}}{8}}+$\par${\frac {3765282218997607935\,
{p}^{7}{h}^{35}}{512}}+{\frac {1477531861770819915\,{p}^{8}{h}^{33}}{
2048}}+{\frac {313794615563382465\,{p}^{9}{h}^{31}}{4096}}+{\frac {
145420138664361261\,{p}^{10}{h}^{29}}{16384}}+$\par${\frac {
148925406295325835\,{p}^{11}{h}^{27}}{131072}}+{\frac {
85668471069217551\,{p}^{12}{h}^{25}}{524288}}+{\frac {
13954513249092609\,{p}^{13}{h}^{23}}{524288}}+{\frac {
10148990043310371\,{p}^{14}{h}^{21}}{2097152}}+$\par${\frac {
31811088804006807\,{p}^{15}{h}^{19}}{33554432}}+{\frac {
25711943256475155\,{p}^{16}{h}^{17}}{134217728}}+{\frac {
10166014410124383\,{p}^{17}{h}^{15}}{268435456}}+{\frac {
7042540335027723\,{p}^{18}{h}^{13}}{1073741824}}+$\par${\frac {
4488743637561879\,{p}^{19}{h}^{11}}{8589934592}}-{\frac {
14409338464147941\,{p}^{20}{h}^{9}}{34359738368}}-{\frac {
30196613342204865\,{p}^{21}{h}^{7}}{68719476736}}-{\frac {
198392485810791933\,{p}^{22}{h}^{5}}{549755813888}}-$\par${\frac {
5911762450029857199\,{p}^{23}{h}^{3}}{17592186044416}}-{\frac {
31891120592767324251\,{p}^{24}h}{70368744177664}}
$}.

\end{thm}
\begin{rem}
Note that the coefficients of the monomial $h^{N-2j}p^j$  of $ A_{N}$
are not well-defined for
$j>(N-5)/4$  because they do not lie in the stable range.
We would  like to thank Dieter Kotschick for pointing this out.
One might however view the above formulas  as a definition
of these addititonal terms.
One also sees that many of the invariants  out of the stable range are
negative whereas all those inside the stable range are positive
(this was also pointed out to us by Dieter Kotschick).
\end{rem}

\subsection{The $SO(3)$-case}
We consider first the $SO(3)$-invariants on $\widehat \P_2$ with respect to
Chern classes
$(H,c_2)$, and we put $N:=4c_2-4$.
For $0<\epsilon <<1$ the polarisation $L_\epsilon:=H-\epsilon E$ of
$\widehat \P_2$
 lies in a chamber of type $(H,c_2)$  related to
 the polarisation $H$ of $\P_2$. Thus (\ref{subl}) gives
$$\Phi^{\P_2,H}_{H,N}(\check H^{N-2r}\pt^r)
=\Phi^{\widehat\P_2,L_{\epsilon}}_{H,N}(\check H^{N-2r}\pt^r).$$
Putting $L_{1-\epsilon}:=H-(1-\epsilon) E$ we obtain
$\Phi^{\widehat \P_2,L_{1-\epsilon}}_{H}(\check H^{N-2r}\pt^r)=0$.
Thus we get
$$\Phi^{\P_2,H}_{H}(\check H^{N-2r} \pt^r)=
\sum_{\xi\in W^{\widehat\P_2}_{H,c_2}(H-E,H)}(-1)^{e_{\xi,N}}
\delta_{\xi,N}(\check H^{N-2r}\pt^r),$$
and, using lemma \ref{wallp1p1}, we can again carry out the computation
with Botts formula.
For $N:=4i$ we denote
$$B_N:=2^{2i}\sum_{j=0}^{2i} \Phi^{\P_2,H}_{H}(\check
H^{N-2j}\pt^j)h^{N-2j}p^j.$$
Then we obtain:

\begin{thm}\label{p2so}
The $SO(3)$-Donaldson invariants of $\P_2$ are
\par
{\tolerance 10000\small\raggedright
$B_{0}=1$
, $B_{4}=3\,{h}^{4}+5\,{h}^{2}p+19\,{p}^{2}$,
$B_{8}/8=29\,{h}^{8}+19\,{h}^{6}p+17\,{h}^{4}{p}^{2}+23\,{h}^{2}{p}^{3}+85\,{p}
^{4}$,
\par
$B_{12}=69525\,{h}^{12}+26907\,{h}^{10}p+12853\,{h}^{8}{p}^{2}+7803\,{h}^{6}{p
}^{3}+6357\,{h}^{4}{p}^{4}+8155\,{h}^{2}{p}^{5}+29557\,{p}^{6}$,
\par
$B_{16}/8=6231285\,{h}^{16}+1659915\,{h}^{14}p+519777\,
{h}^{12}{p}^{2}+194439\,{
h}^{10}{p}^{3}+88701\,{h}^{8}{p}^{4}+51027\,{h}^{6}
{p}^{5}+39753\,{h}^
{4}{p}^{6}+49519\,{h}^{2}{p}^{7}+176837\,{p}^{8}$,
\par
$B_{20}=68081556995\,{h}^{20}+13571675125\,{h}^{18}p+3084569555\,{h}^{16}{p}^{
2}+808382629\,{h}^{14}{p}^{3}+247407779\,{h}^{12}{p}^{4}+89811541\,{h}
^{10}{p}^{5}+39553139\,{h}^{8}{p}^{6}+21987589\,{h}^{6}{p}^{7}+
16652099\,{h}^{4}{p}^{8}+20329653\,{h}^{2}{p}^{9}+71741715\,{p}^{10}$,
\par
$B_{24}/8=19355926872345\,{h}^{24}+3046788353175\,{h}^{22}p+535206161485\,{h}^{
20}{p}^{2}+105824308635\,{h}^{18}{p}^{3}+23774344785\,{h}^{16}{p}^{4}+
6132120911\,{h}^{14}{p}^{5}+1838332965\,{h}^{12}{p}^{6}+651103923\,{h}
^{10}{p}^{7}+279395017\,{h}^{8}{p}^{8}+151590087\,{h}^{6}{p}^{9}+
112496445\,{h}^{4}{p}^{10}+135266955\,{h}^{2}{p}^{11}+472659585\,{p}^{
12}$,
\par
$B_{28}=536625215902182969\,{h}^{28}+69259301021976999\,{h}^{26}p+
9817859613586809\,{h}^{24}{p}^{2}+1538955926660199\,{h}^{22}{p}^{3}+
268722697637049\,{h}^{20}{p}^{4}+52689438785319\,{h}^{18}{p}^{5}+
11702994789369\,{h}^{16}{p}^{6}+2974340336103\,{h}^{14}{p}^{7}+
875889126201\,{h}^{12}{p}^{8}+304140743847\,{h}^{10}{p}^{9}+
127923966585\,{h}^{8}{p}^{10}+68135251815\,{h}^{6}{p}^{11}+49776298425
\,{h}^{4}{p}^{12}+59127015975\,{h}^{2}{p}^{13}+204876497145\,{p}^{14}$,
\par
$B_{32}/8=332465777488176686045\,{h}^{32}+36176961518799287203\,{h}^{30}p+
4270943660527526777\,{h}^{28}{p}^{2}+550013108311246927\,{h}^{26}{p}^{
3}+77722220365607813\,{h}^{24}{p}^{4}+12129004922528395\,{h}^{22}{p}^{
5}+2104879834580993\,{h}^{20}{p}^{6}+409294250544727\,{h}^{18}{p}^{7}+
89934657950957\,{h}^{16}{p}^{8}+22556396083123\,{h}^{14}{p}^{9}+
6542216760905\,{h}^{12}{p}^{10}+2235172850335\,{h}^{10}{p}^{11}+
925169690645\,{h}^{8}{p}^{12}+485534741275\,{h}^{6}{p}^{13}+
350230091345\,{h}^{4}{p}^{14}+411833933095\,{h}^{2}{p}^{15}+
1416634092797\,{p}^{16}$,
\par
$B_{36}=17982292064097834276691197\,{h}^{36}+1685376850354867108198203\,{h}^{
34}p+169728914674713290425549\,{h}^{32}{p}^{2}+18446964561578451602667
\,{h}^{30}{p}^{3}+2174127485943121961373\,{h}^{28}{p}^{4}+
279319741333450241435\,{h}^{26}{p}^{5}+39339602087475090285\,{h}^{24}{
p}^{6}+6111138005878747467\,{h}^{22}{p}^{7}+1054025359144892989\,{h}^{
20}{p}^{8}+203321142108471291\,{h}^{18}{p}^{9}+44233113780975117\,{h}^
{16}{p}^{10}+10964566444466603\,{h}^{14}{p}^{11}+3139014782527197\,{h}
^{12}{p}^{12}+1058019835991643\,{h}^{10}{p}^{13}+432158763674797\,{h}^
{8}{p}^{14}+224042778598923\,{h}^{6}{p}^{15}+159901382125437\,{h}^{4}{
p}^{16}+186411458197691\,{h}^{2}{p}^{17}+637107121682253\,{p}^{18}$,
\par
$B_{40}/8=19983831593150830258093037499
\,{h}^{40}+1640698532032417214980201941\,
{h}^{38}p+143617787626796457582947271\,{h}^{36}{p}^{2}+
13451663520190902994423761\,{h}^{34}{p}^{3}+1353428584063925323593987
\,{h}^{32}{p}^{4}+146907352128976242766365\,{h}^{30}{p}^{5}+
17282999997688436388975\,{h}^{28}{p}^{6}+2214864601846913417145\,{h}^{
26}{p}^{7}+310874334747308389131\,{h}^{24}{p}^{8}+48070219333713236901
\,{h}^{22}{p}^{9}+8241254396581767639\,{h}^{20}{p}^{10}+
1577751227160324321\,{h}^{18}{p}^{11}+340134212696649171\,{h}^{16}{p}^
{12}+83440287229631085\,{h}^{14}{p}^{13}+23620292992955391\,{h}^{12}{p
}^{14}+7869891016663881\,{h}^{10}{p}^{15}+3178622918644059\,{h}^{8}{p}
^{16}+1630875748081269\,{h}^{6}{p}^{17}+1153440155417319\,{h}^{4}{p}^{
18}+1334613223327473\,{h}^{2}{p}^{19}+4535236702668195\,{p}^{20}$,
\par
$B_{44}/8=226901192268190530686926956861797\,{h}^{44}+
16542462134525153318253326085835\,{h}^{42}p+$
\par
$1277706977403778580365852666661\,{h}^{40}{p}^{2}+
104862798979925329727378003659\,{h}^{38}{p}^{3}+$\par$
9174416297780080293761973989\,{h}^{36}{p}^{4}+
858689743856030000767365835\,{h}^{34}{p}^{5}+
86310585758469215596920485\,{h}^{32}{p}^{6}+9355633875773319246298315
\,{h}^{30}{p}^{7}+1098557533992391977544805\,{h}^{28}{p}^{8}+$\par$
140418552503311458801355\,{h}^{26}{p}^{9}+19640467303990766625317\,{h}
^{24}{p}^{10}+3023185118099492260555\,{h}^{22}{p}^{11}+$\par$
515310612119604105701\,{h}^{20}{p}^{12}+97958161753205459659\,{h}^{18}
{p}^{13}+20943663715791766949\,{h}^{16}{p}^{14}+$\par$5090445779293122763\,{
h}^{14}{p}^{15}+1426864216020459365\,{h}^{12}{p}^{16}+
470672723850968779\,{h}^{10}{p}^{17}+188268529044707621\,{h}^{8}{p}^{
18}+95733138877112011\,{h}^{6}{p}^{19}+67173305015551205\,{h}^{4}{p}^{
20}+77210866621686475\,{h}^{2}{p}^{21}+261019726029655205\,{p}^{22}$,\par
$B_{48}/64=401623524463671616144253869033873677\,{h}^{48}+
26294009028509419866433950400817907\,p{h}^{46}+$
\par
$1814139310232228402229320933713849\,{p}^{2}{h}^{44}+
132233743700306798807714195145903\,{p}^{3}{h}^{42}+$
\par
$
10210502184961866655190088128661\,{p}^{4}{h}^{40}+
837650587235973991054920612155\,{p}^{5}{h}^{38}+$
\par
$
73245138148540706205679224225\,{p}^{6}{h}^{36}+
6850202117661264825075213975\,{p}^{7}{h}^{34}+$
\par
$
687815006512629065815415005\,{p}^{8}{h}^{32}
+74447573563889724907246275\,{p}^{9}{h}^{30}
+8724562938113746968261705\,{p}^{10}{h}^{28}+$
\par
$
1112238016349638764497855\,{p}^{11}{h}^{26}
+155030522200663787180517\,{p}^{12}{h}^{24}
+23757432754397656762251\,{p}^{13}{h}^{22}
+$
\par
$
4027259666817871766449\,{p}^{14}{h}^{20}+
760537452815695217703\,{p}^{15}{h}^{18}
+161380994483655259053\,{p}^{16}{h}^{16}
+$
\par
$
38900268404198860691\,{p}^{17}{h}^{14}+10809116358008226777
\,{p}^{18}{h}^{12}+3534337177551658959\,{p}^{19}{h}^{10}+$
\par
$
1401766305125084725\,{p}^{20}{h}^{8}+707191549935960795\,{p}^{21}{h}^{
6}+492754565374149825\,{p}^{22}{h}^{4}+563040363143655095\,{p}^{23}{h}
^{2}+1894476461608956285\,{p}^{24}$.
}
\end{thm}

\begin{conj} For all $n$ there is a nonnegative integer $l(n)$ such that
$B_{4n}/2^{l(n)}$ is a polynomial in $h$ and $p$  all of whose coefficients
are odd
positive integers.
\end{conj}

\begin{rem}
$B_0$ to $B_{16}$ were already computed in \cite{K-L} also using blowup and
wall-crossing formulas, showing that $\P_2$ is not of simple type.
Apart from slightly different conventions their results
agree with ours. Their results and the computations of the $SU(2)$
invariants by various other authors have been quite useful to check the
correctness
of our programs -- and thus of the computations in section 3 --
 in earlier
stages of our work.
The conjecture above could already have been made on the basis of their result.
\end{rem}

\section{Wall-crossing formulas}
In our paper \cite{E-G} we formulated a conjecture about the
shape of the wall-crossing  formula,  compatible with the
conjecture of Kotschick and Morgan \cite{K-M}.
Here we state  a slightly stronger form of the conjecture
which is also supported by the computations in
 \cite{E-G}.

\begin{conj}\label{wallconj}  In the polynomial ring on $H^2(S,\Q)$ we have
$$\delta_{\xi,N,r}=(-4)^{-r}\sum_{k=0}^d
{(N-2r)!\over (N-2r-2d+2k)!(d-k)!}Q_{k}(N,d,r,K_S^2)
L_{\xi/2}^{N-2r-2d+2k}q_S^{d-k},$$
where $Q_{k}(N,d,r,K_S^2)$ is a  polynomial of degree
$k$ in $N,d,r,K_S^2$, which is independent of $S$ and $\xi$.
\end{conj}

We now will show, that, assuming the conjecture, we can compute
several of the $Q_{k}(N,d,r,K^2_S)$.
This computation will also give a check of the conjecture in many
specific cases.

For all $i\ge 0$ we put
$$P_i(N,d,r,K^2):={-2d+2N+2K^2-24r+{\frac{13+3i}{2}}\choose i}
+(3N-288r){-2d+2N+2K^2-24r+{\frac{7+3i}{2}}
\choose i-2}.$$

\begin{prop}\label{wallprop1}
If conjecture \ref{wallconj} is true, then
for $i=0,1,2,3,4$ we can write
$Q_{i}(N,d,r,K^2)=P_i(N,d,r,K^2)+R_{i}(N,d,r,K^2),$
where  $R_{i}(N,d,r,K^2)=0$ for $i<2$ and
\newline
{\raggedright
$ R_2(N,d,r,K^2)={\frac {69}{8}}$,
\newline
$ R_3(N,d,r,K^2)= -13\,d+29\,N+17\,K^2-5452\,r+91$,
\newline
$ R_4(N,d,r,K^2)= {\frac {35\,{d}^{2}}{4}}
-{\frac {99\,dN}{2}}-{\frac {51\,d(K^2)}{2}
}+10802\,dr+{\frac {181\,{N}^{2}}{4}}+{\frac {115\,N(K^2)}{2}}-
12050\,Nr+{\frac {67\,{(K^2)}^{2}}{4}}-10898\,(K^2)\,r+169836\,{
r}^{2}-{\frac {537\,d}{4}}+{\frac {2821\,N}{8}}+{\frac {761\,(K^2)}
{4}}-146495\,r+{\frac {72005}{128}}.$}
\end{prop}
\begin{pf}
We assume conjecture \ref{wallconj}.
Let $X=\P_1\times \P_1$ or a blow up of $\P_1\times \P_1$ in finitely many
points. We denote by $F$ and $G$ the pullbacks of the fibres of the
two projections from $\P_1\times \P_1$ to $\P_1$.
For a class $\xi=F-sG$  in $H^2(X,\Z)$ defining a wall
and an integer $d\ge 0$,
let $N:=4d+2s-3$. Then on $X$ we can determine the coefficients $a_{k}$
of $L_{\xi/2}^{N-2r-2d+2k}q_X^{d-k}$ in $\delta_{\xi,N,r}$ as follows:
We can  assume that $X$ has an action of $\C^*$ with finitely many fixpoints
as in \ref{torS}.
For $x$ an indeterminant  we put $\alpha:=-xF+G$ and compute the polynomial
$\delta_{\xi,N,r}(\alpha^{N-2r})$ in $x$ using the Bott residue formula
on $\Hilb^d(X\sqcup X)$, from which we can compute the
$a_k$.

Now we can compute the polynomials $Q_k(N,d,r,K_X^2)$ as follows:
We consider all nonnegative integers $d,w,b,r$ with $d\ge k$ and $d+w+b+r\le
2k$.
Let $X$ be the blow up of $\P_1\times \P_1$ in $b$ points.
With the method of the last paragraph we compute the coefficient $c_{d,w,b,r}$
of
$L_{\xi/2}^{N-2r-2d+2k}q_X^{d-k}$ in $\delta_{\xi,N,r}$ on $X$,
where $N:=4d+2w+1$ and $\xi=F-(w+2)G$. Using all the  $c_{d,w,b,r}$ we obtain
a system of ${k+4 \choose 4}$ linear equations for the coefficients of
the $N^id^jr^t(K_S^2)^s$ (with $0\le i+j+t+s\le k$)  in $Q_k(N,d,r,K_S^2)$.
Solving this system of equations we obtain our result.
All the computations are again carried out using a suitable Maple program.
\end{pf}

The formulas suggest the following conjecture:

\begin{conj} \label{wallconj1}
\begin{enumerate}
\item For all $i$ the polynomial $Q_i(N,d,r,K^2)$ is of the form
$Q_{i}(N,d,r,K^2)=P_{i}(N,d,r,K^2)+R_{i}(N,d,r,K^2),$
where $R_{i}(N,d,r,K^2)$ is a polynomial in
$N,d,r,K^2$ of degree $i-2$.
\item If we view $R_i(N,d,r,K^2)$ as a polynomial in
$N,-d,-r,K^2$, then all its coefficients are positive and the same is true
for $Q_i(N,d,r,K^2)$.
\end{enumerate}
\end{conj}

\begin{prop}\label{wallprop2}
If conjecture \ref{wallconj} and part (1) of conjecture \ref{wallconj1}
are true then
\newline
{\raggedright
$R_5(N,d,r,K^2)=
11\,{{(K^2)}}^{3}+57\,{{(K^2)}}^{2}N-25\,{{(K^2)}}^{2}d-10892\,{{
(K^2)}}^{2}r+90\,{(K^2)}\,{N}^{2}-98\,{(K^2)}\,Nd-24088\,{(K^2)}\,
Nr+17\,{(K^2)}\,{d}^{2}+21592\,{(K^2)}\,dr+339600\,{(K^2)}\,{r}^{2}
+44\,{N}^{3}-82\,{N}^{2}d-13304\,{N}^{2}r+41\,N{d}^{2}+23896\,Ndr+
363792\,N{r}^{2}-3\,{d}^{3}-10700\,{d}^{2}r-338448\,d{r}^{2}-2525760\,
{r}^{3}+198\,{{(K^2)}}^{2}+744\,{(K^2)}\,N-276\,{(K^2)}\,d-303600\,
{(K^2)}\,r+618\,{N}^{2}-612\,Nd-333888\,Nr+78\,{d}^{2}+301008\,dr+
5457888\,{r}^{2}+1213\,{(K^2)}+2506\,N-729\,d-3101884\,r+2490
$,
\newline
\smallskip
$R_6(N,d,r,K^2)=   {\frac {65\,{{(K^2)}}^{4}}{12}}+{\frac
{113\,{{(K^2)}}^{3}N}{3}}-{
\frac {49\,{{(K^2)}}^{3}d}{3}}-{\frac {21772\,{{(K^2)}}^{3}r}{3}}+{
\frac {179\,{{(K^2)}}^{2}{N}^{2}}{2}}-97\,{{(K^2)}}^{2}Nd-24076\,{{
(K^2)}}^{2}Nr+{\frac {33\,{{(K^2)}}^{2}{d}^{2}}{2}}+21580\,{{(K^2)}
}^{2}dr+339528\,{{(K^2)}}^{2}{r}^{2}+{\frac {263\,{(K^2)}\,{N}^{3}}{
3}}-163\,{(K^2)}\,{N}^{2}d-$
\par
$26596\,{(K^2)}\,{N}^{2}r+81\,{(K^2)}\,N{
d}^{2}+47768\,{(K^2)}\,Ndr+727440\,{(K^2)}\,N{r}^{2}-{\frac {17\,{
(K^2)}\,{d}^{3}}{3}}-21388\,{(K^2)}\,{d}^{2}r-676752\,{(K^2)}\,d{r}
^{2}-5050944\,{(K^2)}\,{r}^{3}+{\frac {365\,{N}^{4}}{12}}-{\frac {247
\,{N}^{3}d}{3}}-{\frac {29332\,{N}^{3}r}{3}}+{\frac {147\,{N}^{2}{d}^{
2}}{2}}+26404\,{N}^{2}dr+389208\,{N}^{2}{r}^{2}-{\frac {65\,N{d}^{3}}{
3}}-23692\,N{d}^{2}r-725136\,Nd{r}^{2}-5327424\,N{r}^{3}+{\frac {{d}^{
4}}{12}}+{\frac {21196\,{d}^{3}r}{3}}+337224\,{d}^{2}{r}^{2}+5041728\,
d{r}^{3}+24147648\,{r}^{4}+{\frac {821\,{{(K^2)}}^{3}}{6}}+{\frac {
3127\,{{(K^2)}}^{2}N}{4}}-{\frac {565\,{{(K^2)}}^{2}d}{2}}-314198\,{
{(K^2)}}^{2}r+1306\,{(K^2)}\,{N}^{2}-{\frac {2567\,{(K^2)}\,Nd}{2}}
-691258\,{(K^2)}\,Nr+{\frac {309\,{(K^2)}\,{d}^{2}}{2}}+623020\,{
(K^2)}\,dr+11250984\,{(K^2)}\,{r}^{2}+{\frac {7987\,{N}^{3}}{12}}-
1154\,{N}^{2}d-380192\,{N}^{2}r+{\frac {2007\,N{d}^{2}}{4}}+685594\,Nd
r+12076284\,N{r}^{2}-{\frac {53\,{d}^{3}}{6}}-308822\,{d}^{2}r-
11204904\,d{r}^{2}-93985056\,{r}^{3}+{\frac {249499\,{{(K^2)}}^{2}}{
192}}+{\frac {529411\,{(K^2)}\,N}{96}}-{\frac {48265\,{(K^2)}\,d}{32
}}-{\frac {156136081\,{(K^2)}\,r}{24}}+{\frac {944107\,{N}^{2}}{192}}
-{\frac {134657\,Nd}{32}}-{\frac {170651209\,Nr}{24}}+{\frac {41627\,{
d}^{2}}{192}}+{\frac {154931089\,dr}{24}}+{\frac {538005841\,{r}^{2}}{
4}}+{\frac {1027343\,{(K^2)}}{192}}+{\frac {5033035\,N}{384}}-{\frac
{406799\,d}{192}}-{\frac {1032193219\,r}{16}}+{\frac {7872921}{1024}}.$}
\end{prop}

\begin{pf}
The method of the proof is very similar to that of proposition
\ref{wallprop1}.
Using conjecture \ref{wallconj1} we can reduce the number of equations
of the linear system we have to solve. For the computation of $R_k$
we have in the notations of
the proof of proposition \ref{wallprop1} only to consider
$d,w,b,r$ with $d\ge k$ and $d+w+b+r\le
2k-2$.
\end{pf}

\begin{prop}
Assume  conjecture \ref{wallconj}. Then for   $d\le 8$ and  $K^2=8$
also conjecture \ref{wallconj1} holds.
Furthermore we have
\newline
{\raggedright
$R_7(N,7,r,8)={\frac {242\,{N}^{5}}{15}}-{\frac
{16136\,{N}^{4}r}{3}}+277280\,{N}^{3
}{r}^{2}-5613696\,{N}^{2}{r}^{3}+50448384\,N{r}^{4}-{\frac {843153408
\,{r}^{5}}{5}}+584\,{N}^{4}-{\frac {926104\,{N}^{3}r}{3}}+14133648\,{N
}^{2}{r}^{2}-215219520\,N{r}^{3}+1068715008\,{r}^{4}+{\frac {49259\,{N
}^{3}}{6}}-{\frac {27036430\,{N}^{2}r}{3}}+327411912\,N{r}^{2}-
2834114784\,{r}^{3}+55637\,{N}^{2}-{\frac {491711972\,Nr}{3}}+
3453170368\,{r}^{2}+{\frac {5432731\,N}{30}}-{\frac {7607891522\,r}{5}
}+226878$,
\newline
$R_7(N,8,r,8)={\frac {242\,{N}^{5}}{15}}-{\frac
{16136\,{N}^{4}r}{3}}+277280\,{N}^{3
}{r}^{2}-5613696\,{N}^{2}{r}^{3}+50448384\,N{r}^{4}-{\frac {843153408
\,{r}^{5}}{5}}+526\,{N}^{4}-{\frac {867848\,{N}^{3}r}{3}}+13357680\,{N
}^{2}{r}^{2}-204584256\,N{r}^{3}+1020478464\,{r}^{4}+{\frac {40547\,{N
}^{3}}{6}}-{\frac {24608014\,{N}^{2}r}{3}}+301879416\,N{r}^{2}-
2636445408\,{r}^{3}+{\frac {85365\,{N}^{2}}{2}}-{\frac {444885076\,Nr}
{3}}+3161672456\,{r}^{2}+{\frac {1989263\,N}{15}}-{\frac {6895690692\,
r}{5}}+162324$,
\newline
$R_8(N,8,r,8)={\frac {69\,{N}^{6}}{10}}-{\frac
{35428\,{N}^{5}r}{15}}+147976\,{N}^{4
}{r}^{2}-3940224\,{N}^{3}{r}^{3}+52663680\,{N}^{2}{r}^{4}-{\frac {
1749924864\,N{r}^{5}}{5}}+{\frac {4618156032\,{r}^{6}}{5}}+{\frac {
17993\,{N}^{5}}{60}}-164166\,{N}^{4}r+9819288\,{N}^{3}{r}^{2}-
222171264\,{N}^{2}{r}^{3}+2193801408\,N{r}^{4}-{\frac {39855794688\,{r
}^{5}}{5}}+{\frac {342159\,{N}^{4}}{64}}-{\frac {75235799\,{N}^{3}r}{
12}}+{\frac {675173309\,{N}^{2}{r}^{2}}{2}}-5810679060\,N{r}^{3}+
32006797596\,{r}^{4}+{\frac {9540439\,{N}^{3}}{192}}-{\frac {678142435
\,{N}^{2}r}{4}}+{\frac {28388710281\,N{r}^{2}}{4}}-68700098862\,{r}^{3
}+{\frac {646487951\,{N}^{2}}{2560}}-{\frac {602240211743\,Nr}{192}}+{
\frac {11796298005871\,{r}^{2}}{160}}+{\frac {10191068747\,N}{15360}}-
{\frac {19195182347591\,r}{640}}+{\frac {23061793325}{32768}}$.}
\end{prop}
\begin{pf}
The method is again similar to that of the proof of
proposition \ref{wallprop1}.
Now we carry out our computations on $X=\P_1\times\P_1$.
In the notation of the proof of proposition \ref{wallprop1}
we have therefore $b=0$. For the computation of $R_k$ we
consider nonnegative integers $d,w,r$ with
$d+w+r\le 2k$ and $k\le d \le 8$.
\end{pf}

\begin{prop}\label{wallrul}
Let $S$ be a rational ruled surface, then for $N\le 40$ and $d\le 8$
the conjectures \ref{wallconj} and \ref{wallconj1} are correct
(and therefore also all the formulas above).
\end{prop}

\begin{pf}
Any rational ruled surface $X$ is a degeneration of either
or $\P_1\times \P_1$ or $\widehat \P_2$, and under the
 degeneration the ample cone of $X$
corresponds
to a subcone of the ample cone of $\P_1\times \P_1$
(resp. $\widehat \P_2$).
Therefore it is enough to prove the result for $\P_1\times \P_1$
and $\widehat \P_2$.
We let  $c_1$ run through $0,F,G,F+G$ on $\P_1\times \P_1$  and through
$0,H,E,F$ on $\widehat \P_2$ ($F=H-E$).
For $S=\P_1\times\P_1$ and $S=\widehat \P_2$ we consider  for all
integers $d$ with $0\le d\le 8$
the set $W_{S,d}$
of all  classes $\xi$, which define a wall of type
$(c_1,c_2)$,  such that $N:=4d-\xi^2-3\le 40$ and $\<\xi\cdot F\><0$.
It is easy to see that
\begin{eqnarray*}
W_{\P_1\times\P_1,d}&=&\Big\{\xi=aF-bG\Bigm| a>0,b>0,2ab\le 40-4d+3\Big\},\\
W_{{\widehat\P_2},d}&=&\Big\{\xi=bH-aE\Bigm| a>b>0, a^2-b^2\le 40 -4d+3\}.
\end{eqnarray*}
For all $d\le 8 $ and all $\xi\in W_{S,d}$ we again compute
all the coefficients of $L_{\xi/2}^{N-2r-2d+2k}q_{S}^{d-k}$
with the method of the first paragraph of the proof of proposition
\ref{wallprop1}.
\end{pf}

\section{The Donaldson Invariants of birationally ruled surfaces.}
In this section we will show
that our algorithm for computing the wall-crossing formula $\delta_{\xi,N}$
and the blowup formulas enable us to compute all the
Donaldson invariants of all rational surfaces $X$ for all generic polarisations
lying in a suitable subcone of the
 ample cone of $X$.

\subsection{The case of rational ruled surfaces}
In this case we can
indeed determine the Donaldson invariants for all generic polarisations.
For simplicity we will only compute the restriction of the Donaldson invariants
to $\Sym^N(H_2(S,\Q))$.
In \cite{K-L} some invariants of $\P_1\times \P_1$ have been computed also
using
blowup and wall-crossing formulas. Their results show e.g.
that there is no chamber, for which
$\P_1\times \P_1$
is of simple type. Our results again agree with theirs
and earlier results e.g. in \cite{L-Q}.

\begin{thm}\label{donrul}
Let $S$ be a rational ruled surface, $F$ the class of a fibre
and $q_S$ the quadratic form on $H_2(S,\Z)$.
We denote by $F_\epsilon$ the polarisation $F+\epsilon E$, where
$E$ is the class of a section with nonpositive selfintersection.

\noindent $(1)$
For $\epsilon>0$ sufficiently small
we have $\Phi^{S,F_\epsilon}_{E,N}=0$
and
$\Phi^{S,F_\epsilon}_{E+F,N}=0.$

\noindent $(2)$
For $\epsilon>0$ sufficiently small
we have for $E_N:=\Phi^{S,F_\epsilon}_{0,N,0}$:
\newline
{\raggedright
$E_5=-L_F^5+5/2L_F^3q_S-5/2L_Fq_S^2$,
\newline
$E_9=40L_F^9-108L_F^7q_S+108L_F^5q_S^2-42L_F^3q_S^3$,
\newline
$E_{13}=-9345L_F^{13}+26949L_F^{11}q_S-31590L_F^9q_S^2
+18018L_F^7q_S^3-4290L_F^5q_S^4$,
\newline
{\small
$E_{17}=7369656L_F^{17}-22136040L_F^{15}q_S
+28474320L_F^{13}q_S^2-19734960L_F^{11}q_S^3
+7425600L_F^9q_S^
4-1225224L_F^7q_S^5$,
\newline
$E_{21}=-14772820744L_F^{21}+45586042992L_F^{19}q_S-62181472500L_F^{17}q_S^2
+48231175860L_F^{15}q_S^3-$
\newline
$
22562971200L_F^{13}q_S^4+6074420688L_F^{11}q_S^5-740703600L_F^9q_S^6$,
\newline
$E_{25}=63124363433664L_F^{25}-198545836440000L_F^{23}q_S
+281925714232800L_F^{21}q_S^2-
235199340734400L_F^{19}q_S^3+$
\newline
$
125056219068000L_F^{17}q_S^4-42588214875360L_F^{15}q_S^5+
8649138960000L_F^{13}q_S^6-813136737600L_F^{11}q_S^7$,
\newline
$E_{29}=-509894102555251905L_F^{29}+1626742370158553130L_F^{27}q_S
-2378321090933081112L_F^{25}q_S^2+$
\newline
$
2087846466793743600L_F^{23}q_S^3-1207966082767844400L_F^{21}q_S^4
+473530658232013200
L_F^{19}q_S^5-$
\newline
$
123363365393268000L_F^{17}q_S^6+19623703009790880L_F^{15}q_S^7-
1467326424564000L_F^{13}q_S^8,$\newline
$E_{33}=7135482220088837442520\,L_F^{33}-23016295766978863295760\,L_F^{31}q_S+
34404291587748659734080\,L_F^{29}q_S^{2}-$
\newline
$31360607908598315276160\,L_F^{27}q_S^{3}+19266231547036209415680\,L_F^{25}
q_S^{4}-
8299005150626510918400\,L_F^{23}q_S^{5}+$
\newline
$
2515398487826672448000\,L_F^{
21}q_S^{6}-519339581441771650560\,L_F^{19}q_S^{7}+66567414222758592000
\,L_F^{17}q_S^{8}-$
\newline
$
4055565288690115200\,L_F^{15}q_S^{9}.$}}

\noindent $(3)$ For $\epsilon>0$ sufficiently small and all
$N\le 33$ we have, writing $\Phi^{S,F_\epsilon}_{0,N,0}$
 as a polynomial
$E_N(L_F,q_S)$ in $L_F$ and $q_S$,
$$\Phi^{S,F_\epsilon}_{F,N,0}=E_N(L_F,q_S)
-E_N(L_F/2,q_S).$$
\end{thm}
\begin{pf}
(1)  is just lemma \ref{vancham}.

\noindent (2) and (3): We fix $N:=4c_2-3$ with $c_2>1$.
We will just compute the corresponding Donaldson invariants explicitely.
As any  pair $(S,L)$ consisting of a Hirzebruch surface $S=\Sigma_n$ and
$L=aF+bE\in Pic(S)$ (where $E$ is a section with selfintersection $-n\le 0$)
can be deformed to either
 $(\P_1\times \P_1, aF+b(E-nF/2))$ or $(\widehat \P_2,aF+b(E-(n-1)F/2))$
we see that we can assume that $S=\P_1\times \P_1$ or $S=\widehat \P_2$
and $c_1=F$ or $c_1=0$.

(a) $S=\P_1\times\P_1$, $c_1=F$ (we will always denote by $F$ and $G$
the fibres of the projections to the two factors).
By (1) we have for $\epsilon>0$ sufficiently small
$\Phi^{\P_1\times\P_1,G+\epsilon F}_{F,N}=0$.
Therefore we get
$$\Phi^{\P_1\times\P_1,F_{\epsilon}}_{F,N}=
-\sum_{\xi\in W^{\P_1\times \P_1}_{F,c_2}(F,G)}
(-1)^{e_{\xi,N}}\delta_{\xi,N}.$$
So the invariants can be computed using proposition \ref{wallrul} and
remark \ref{wallp1p1}.

(b) $S=\widehat \P_2$, $c_1=F$.
Let $\epsilon>0$ be sufficiently small.
By the $SO(3)$-blow up formula
we have for all $\alpha\in A_{N-i}(\P_2)$:
$$\Phi^{\widehat \P_2,H-\epsilon E}_{F,N}(\check E^i\alpha)
=\Phi^{\P_2,H}_H(S_i(\pt)\alpha),$$
and the  $SO(3)$-invariants of $\P_2$ have been determined in theorem
\ref{p2so}.
Therefore
$$\Phi^{\widehat \P_2,F_\epsilon
}_{F,N}(\check E^i\alpha)=\Phi^{\P_2,H}_H(S_i(\pt)\alpha)
-\sum_{\xi\in  W^{\widehat \P_2}_{F,c_2}(F,H)}
(-1)^{e_{\xi,N}}\delta_{\xi,N}(\check E^i\alpha),$$
So the sum can be computed using proposition \ref{wallrul}  and remark
\ref{wallp1p1}.

(c) $S=\widehat\P_2$, $c_1=0$.
Let $\widetilde \P_2$ be the blowup of $\P_2$ in two points
with exceptional divisors $E_1$ and $E_2$. Then $\widetilde \P_2$ is also
the blow up of $\P_1\times \P_1$ in a point. We denote the exceptional divisor
by $E$. We denote by $F$ the pullback of $F=H-E_1$ from $\widehat \P_2$
(which coincides with the pullback of $F$ from $\P_1\times \P_1$).
We have $F=E_2+E$. We also denote by $G$ the pullback of $G$ from
$\P_1\times \P_1$ and have $G=E_1+E$. For $1>>\epsilon>>\mu>0$,
let $H_1:=F+\epsilon G-\mu E$ and $H_2:=F+\epsilon G- (\epsilon-\mu)E$.
Then $H_2$ is a polarisation of $\widetilde \P_2$
which lies in a $(E_2,c_2)$-chamber related to the
$(0,c_2)$-chamber of $F+\epsilon E_1$ on $\widehat \P_2$.
Thus by the $SO(3)$-blowup formula we have
$$\Phi^{\widehat \P_2,F+\epsilon E_1}_{0,N}(\check F^i\check E_1^{N-i})=
-\Phi^{\widetilde\P_2,H_2}_{E_2,N+1}(\check E_2 \check F^i (\check G-\check
E)^{N-i})=
-\Phi^{\widetilde\P_2,H_2}_{F-E,N+1}
((\check F-\check E) \check F^i (\check G-\check E)^{N-i}).$$
We have
$$
\Phi^{\widetilde\P_2,H_1}_{F-E,N+1}- \Phi^{\widetilde\P_2,H_2}_{F-E,N+1}
=\sum_{\xi\in W^{\widetilde\P_2}_{F-E,c_2}(H_2,H_1)}
(-1)^{e_{\xi,N+1}}\delta_{\xi,N+1},$$ and
 for $\epsilon$ sufficiently small is is easy to see that
$$W^{\widetilde\P_2}_{F-E,c_2}(H_2,H_1)=
\big\{(2a-1)F-(2b-1)E\bigm| b>a>0, b(b-1)\le c_2\big \}.$$
So
$\Phi^{\widetilde\P_2,H_1}_{F-E,N+1}- \Phi^{\widetilde\P_2,H_2}_{F-E,N+1}$
can be computed by the Bott residue formula.
Finally  $H_1$ lies in a $(F-E,c_2)$-chamber on $\widetilde \P_2$
 related to the
$(F,c_2)$-chamber of $F+\epsilon G$ on $\P_1\times  \P_1$.
So, by the $SO(3)$-blowup formula (with exceptional divisor $E$),
we get for $\alpha\in A_{N+1-i}( \P_1\times \P_1)$
$$\Phi^{\widetilde\P_2,H_1}_{F-E,N+1}(\alpha \check E^i)=
\Phi^{ \P_1\times\P_1,F+\epsilon E}_{F}(\alpha S_i(\pt)),$$ and the
last
is computed by the method of (a).
Now we put everything together to get our result.

(d) $S=\P_1\times \P_1$, $c_1=0$.
This case is very similar to (c), only with the role of $\P_1\times\P_1$ and
$\widehat \P_2$ exchanged. We use the same notations as in (c).
Now  $H_1$ is a polarisation of $\widetilde \P_2$
which lies in a $(E,c_2)$-chamber related to the
$(0,c_2)$-chamber of $F+\epsilon G$ on $\P_1\times \P_1$.
Thus by the $SO(3)$-blowup formula we have
$$\Phi^{\P_1\times\P_1,F+\epsilon G}_{0,N}(\check F^i\check G^{N-i})=
-\Phi^{\widetilde\P_2,H_1}_{E,N+1}(\check E \check F^i \check G^{N-i})=
-\Phi^{\widetilde\P_2,H_1}_{F-E_2,N+1}
((\check F-\check E_2) \check F^i (\check F+\check E_1-\check E_2)^{N-i}).$$
We have
$$
\Phi^{\widetilde\P_2,H_1}_{E,N+1}- \Phi^{\widetilde\P_2,H_2}_{E,N+1}
=\sum_{\xi\in W^{\widetilde\P_2}_{E,c_2}(H_2,H_1)}
(-1)^{e_{\xi,N+1}}\delta_{\xi,N+1},$$ and
 for $\epsilon$ sufficiently small is is easy to see that
$$W^{\widetilde\P_2}_{E,c_2}(H_2,H_1)=
\big\{(2aF-(2b-1)E\bigm| b-1/2>a>0, b(b-1)\le c_2\big \}.$$
So
$\Phi^{\widetilde\P_2,H_1}_{E,N+1}- \Phi^{\widetilde\P_2,H_2}_{E,N+1}$
can be computed by the Bott residue formula.
Finally  $H_2$ lies in a $(E,c_2)$-chamber on $\widetilde \P_2$
 related to the
$(F,c_2)$-chamber of $F+\epsilon E_1$ on $\widehat \P_2$.
So, by the $SO(3)$-blowup formula (with exceptional divisor $E_2$),
we get for $\alpha\in A_{N+1-i}(\widehat \P_2)$
$$\Phi^{\widetilde\P_2,H_2}_{E,N+1}(\alpha \check E_2^i)=
\Phi^{\widehat \P_2,F+\epsilon E}_{F}(\alpha S_i(\pt)),$$ and the last
is computed by the method of (b).
\end{pf}

\begin{conj} For $S$ a rational ruled surface we have in the
notation of theorem \ref{donrul} for all $N=4c_2 -3 $ with $c_2\ge 2$:
\begin{enumerate}
\item $\Phi^{S,F_\epsilon}_{0,N,0}$ and $\Phi^{S,F_\epsilon}_{F,N,0}$
are polynomials $E_{0,N}(L_F,q_S)$ and $E_{F,N}(L_F,q_S)$ in $L_F$ and $q_S$,
which are
independent of $S$.
\item
$E_{0,N}(L_F,q_S)$ and $E_{F,N}(L_F,q_S)$
are divisible by $L_F^{N-2c_2}$.
\item
$E_{F,N}(L_F,q_S)=E_{0,N}(L_F,q_S)
-E_{0,N}(L_F/2,q_S).$
\end{enumerate}
\end{conj}

\begin{rem}
We keep the notation of theorem \ref{donrul}.
Notice that theorem \ref{donrul} and proposition \ref{wallrul}
 determines
all the $SU(2)$- and $SO(3)$- Donaldson invariants of
a rational ruled surface
$S$ of degree at most $35$ for all generic polarisations:
Fix $(c_1,c_2)$ and put $N:=4c_2-c_1^2-3$.
If  $L$ is a generic polarisation then
$$\Phi^{L,S}_{c_1,N}=\Phi^{F_\epsilon,S}_{c_1,N}+\sum_{\xi\in
W^S_{c_1,c_2}(F,L)}
(-1)^{e_{\xi,N}}\delta_{\xi,N}.$$
This sum is given
for $N\le 35$ by
theorem \ref{donrul}, remark \ref{wallp1p1}  and proposition \ref{wallrul}.
\end{rem}

\medskip

\subsection{The Donaldson invariants of blowups of $\P_2$}

We want to finish by showing that our methods give an
algorithm for computing all the Donaldson invariants for all rational
surfaces $X$ at least for polarisations lying in a reasonably big subcone
$C^{g}$
of the  ample cone $C_X$ of $X$.
In \cite{K-L} it is shown that
the Donaldson invariants of $\P_2$ and $\P_1\times \P_1$ can be determined
from the wall-crossing formulas on some blowups, and our results can be seen as
a generalization of this.

A rational surface $X$, which is neither $\P_2$ nor ruled can be deformed
to a a blowup $\P_2(x_1,\ldots x_r)$ of $\P_2$ in finitely many general points.
Under this deformation  $C_X$ corresponds to a
(in general strict) subcone of the ample cone
$C_{\P_2(x_1,\ldots x_r)}$. We can therefore restrict our attention
to  $X=\P_2(x_1,\ldots x_r)$.

\begin{thm}
There exists an algorithm computing all the $SU(2)$- and $SO(3)$-Donaldson
invariants of $\P_2(x_1,\ldots x_r)$ with respect to all generic polarisations
in a nonempty open subcone $C^{g}$ of the ample cone of $\P_2(x_1,\ldots x_r)$.
\end{thm}

\begin{pf}
Let $S=Y_r$, where $Y_0=\P_2$ and $Y_i$ is obtained from $Y_{i-1}$
by blowing up a point such  that each $Y_i$ carries an action
of an algebraic $2$-torus with finitely many fixpoints
satisfying conditions (1) and (2) of section
\ref{torS}. This
just means that each $Y_i$ is obtained from $Y_{i-1}$ by blowing up
a fixpoint.
We can deform $S$ to
$\P_2(x_1,\ldots x_r)$, but under this deformation the
good ample cone  $C_S^g$ of $S$ will in general correspond to a proper subcone
$C^g$ of $C_{\P_2(x_1,\ldots x_r)}$.
Note that $C^g$  contains
a neighbourhood of the hyperplane class $H$.
It is enough to prove that there is such an algorithm computing  the Donaldson
invariants of $S$ for all generic polarisations in $C_S^g$.

Fix $c_1\in Pic(S)$ and $c_2\in H^2(S,\Z)$. Let $N:=4c_2-c_1^2-3$.
Let  $H_1$ and $H_2$ be two good generic polarisations  of
$S$. Then  by lemma \ref{finwall} the set
$W_{c_1,c_2}(H_1,H_2)$ is finite and consists only of good walls.
Therefore
$$\Phi_{c_1,N}^{S,H_2}=\Phi_{c_1,N}^{S,H_1}+\sum_{\xi\in
W^S_{c_1,c_2}(H_1,H_2)}
(-1)^{e_{\xi,N}}\delta_{\xi,N},$$
and all the $\delta_{\xi,N}$ can be determined explicitely by applying
the Bott residue formula.
So it is enough to determine $\Phi_{c_1,N}^{S,H_0}$ for one good polarisation
$H_0$. We will denote by $E_1,\ldots,E_r$ the exceptional divisors of
$S$ over $\P_2$.

{\it First case $c_1\ne 0$.}
Denote $c_1=a H+b_1 E_1+\ldots + b_r E_r$, with each of
$a,b_1,\ldots,b_r$ lying in $\{ 0,1\}$. We denote
$D_i:=a H+b_1 E_1+\ldots b_i + E_i$

By reordering the
$E_i$ we can assume $a\ne 0$ or $b_1\ne 0$.
Let $F:=H-E_1$. Then for $1>>\epsilon>>\delta_2>>\ldots >>\delta_r$
the divisors
$H_{i}:=F+\epsilon E_1-(\delta_2 E_2+\ldots+\delta_i E_i) $
are polarisations on $Y_i$
lying  in a chamber of
type $(D_i,c_2)$ related to the chamber of type $(D_{i-1},c_2)$ of
$H_{i-1}:=F+\epsilon E_1-(\delta_2 E_2+\ldots\delta_{i-1} E_{i-1})$
on $Y_{i-1}$.
So the blowup formulas give
$\Phi^{Y_i,H_{i}}_{D_i}(\alpha \check E_i^j)
=\Phi^{Y_{i-1},H_{i-1}}_{D_{i-1}}(\alpha S_j(\pt))$
if $b_i=1$ (resp. $\Phi^{Y_i,H_{i}}_{D_i}(\alpha \check
E_i^j)=\Phi^{Y_{i-1},H_{i-1}}_{D_{i-1}}(\alpha B_j(\pt))$ if $b_i=0$)
for all $\alpha\in A_{N-j}(Y_{i-1})$.
The proof of theorem \ref{donrul} gives an algorithm for computing
$\Phi^{Y_1,H_{1}}_{D_1}(\alpha)$ for all $\alpha\in A_*(Y_{1})$.
Thus by induction we get the desired algorithm.

{\it Second case $c_1=0$.}
Let $\widehat S$ be the blow up of $S$ in a point and let  $E$ be the
exceptional
divisor.
By the first case we have an algorithm for computing the Donaldson invariant
$\Phi_{E,N}^{\widehat S,H_\epsilon}$ for a polarisation
$H_\epsilon=H_0-\epsilon
E$ on $\widehat S$ lying
 in a related chamber to that of a generic good polarisation $H_0$ of $S$.
Then the $SO(3)$-blowup formula gives
$\Phi_{0}^{S,H_0}(\alpha)=\Phi_{E}^{\widehat S,H_\epsilon}(\check E\alpha)$,
and the result follows.
\end{pf}

\end{document}